\documentclass{amsart}

\usepackage{amscd}
\usepackage{amsmath,amssymb,amsfonts}
\usepackage{mathtools}
\usepackage[all, 2cell]{xy}
\usepackage{tikz}
\usepackage{tikz-cd}
\usetikzlibrary{matrix}
\usetikzlibrary{arrows}
\usetikzlibrary{positioning}
\usetikzlibrary{decorations.pathmorphing}

\newcommand{\bC}{\mathbb{C}}
\newcommand{\bZ}{\mathbb{Z}}

\newcommand{\mGL}{{\rm GL}}
\newcommand{\mSL}{{\rm SL}}
\newcommand{\mPGL}{{\rm PGL}}

\newcommand{\mSO}{{\rm SO}}
\newcommand{\mSp}{{\rm Sp}}
\newcommand{\mO}{{\rm O}}

\newcommand{\mGSp}{{\rm GSp}}

\theoremstyle{plain}
\newtheorem{main}{Theorem}

\newtheorem{mainc}[main]{Corollary}

\newtheorem{theorem}{Theorem}[section]
\newtheorem{lemma}[theorem]{Lemma}
\newtheorem{corollary}[theorem]{Corollary}
\newtheorem{proposition}[theorem]{Proposition}

\SelectTips{eu}{12}
\newcommand{\tabincell}[2]{\begin{tabular}{@{}#1@{}}#2\end{tabular}}

\newcommand{\ssmedskip}{\hspace{2ex} \medskip}

\begin{document}

\title{On local and global conjugacy}
% This note is designed for the series talk in USTC 2012 visit, and is to publish in some proper place.

% Reader: Students. All complete details will be included. Useful for LP41.

% The article is a publishable mini note.

\author{Song Wang}
\address{Song Wang, Academy of Mathematics
 and Systematics Sciences, the Morningside Center of Mathematics, 
 the Key Laboratory of Hua, Chinese Academy of Science.}
\thanks{The author was supported in part by One Hundred Talents
Program at CAS, and national 973 project 2013CB834202.}
\email{\texttt{songw1973@amss.ac.cn}}

%    General info
% \subjclass{Primary 22E99, 11F99}
\date{}

% \dedicatory{This paper is dedicated to our advisors.}

% \keywords{Number Theory, Automorphic Forms, Representation Theory}

\maketitle

\tableofcontents

\section{Introduction} \label{S:1}

% Literature

Let $G$ be a reductive linear algebraic group over $\bC$ and $H$ a closed subgroup
of $G$, one may ask to what extent from the representation theory
we can determine $H$. For example, the dimension data for $H$ in $G$
consists of integers $m_{H} (\rho)$ where $\rho$ runs through all finite dimensional
representations of $G$, where $m_{H} (\rho)$ is the multiplicity
of the trivial representation in the restriction of $\rho$ to $H$.
R. ~Langlands (\cite{Langlands-letter2001}, \cite{Langlands2002})
asked, is the dimension data determines the isomorphism class and the conjugacy
class of $H$. In this note, we discuss and classify a special family of counter examples (called LFMO-special
representation, will be defined explicitly later) for $G = \mSO (2 N)$ and $H$ connected
reductive, and give negative answer to Langlands' question in some sense.
In fact, our first family of such examples (Theorem \ref{T:314}, also see
\cite{Wang2007} and \cite{Wang2011}) gave first connected instances of locally conjugate subgroups
of $G = \mSO (2 N)$ failing to be conjugate. Moreover, with sufficient functoriality,
such counter examples will give failure of multiplicity one, and thus got some
attention in the study of beyond endoscopy which is much hotter after the
establishment of the fundamental lemma by NGO (\cite{FLN2010}).

\medskip

In 1990, M. ~Larsen and R. ~Pink studied
the case $G = GL (n)$, and got some results on dimension data (\cite{L-P90}). In fact, they proved,
for general $G$ and connected reductive $H$,
the isomorphism class of $H$ is determined by the dimension data, and when $G = \mGL (n)$
and $H$ embeds into $G$ in an essential way, i.e., irreducible as $n$-dimensional representation,
the conjugacy class of $H$ is determined by the dimension data. \footnote{But for general $H$ in $\mGL (n)$,
the conjugacy of $H$ is not completely determined by the dimension data.
The examples in \cite{L-P90} is not correct,
and the correct example is shown in \cite{A-Y-Y2010}.} For general $G$, this is not the case.
In fact, earlier works exhibited various counter examples
(\cite{Blasius94}, \cite{Larsen94}, \cite{Larsen96},  \cite{Wang2007},
\cite{A-Y-Y2010}, \cite{Yu}, \cite{Wang2011}).

\medskip

% Local VS Conjugacy

The discrepancy between the representation feature and the conjugacy occurs due to various reason,
and one is the ``local-global issue''. Let $G$ as above and $H$, $H'$
two closed subgroups of $G$. We say that $H$ and $H'$ be \emph{locally conjugate}
or \emph{element-wise conjugate} if there is an isomorphism $i: H \overset{\cong}{\longrightarrow} H'$
such that for $h \in H$ in a (Zariski) dense subset of $H$, $h$ and $i (h)$ are conjugate.
We say that $H$ and $H'$ are \emph{globally conjugate} if they are conjugate in $G$.
Also, there are definitions in term of group homomorphisms. Let $\rho, \rho'$: $H \to G$
be two homomorphisms of linear algebraic groups. We say that $\rho$ and $\rho'$
are \emph{locally conjugate} if for $h \in H$ in a (Zariski) dense subset of $H$
$\rho (h)$ and $\rho' (h)$ are conjugate in $G$. We say that $\rho$ and $\rho'$ are
\emph{globally conjugate} if they are conjugate, namely, there exists a $g \in G$
such that $\rho' (h) = g \rho (h) g^{-1}$ for all $h$ in $H$.
We say that $\rho$ and $\rho'$ are \emph{globally conjugate in image}
if $\rho (H)$ and $\rho' (H)$ are conjugate in $G$.

\medskip

Of course, global conjugacy implies local conjugacy. Moreover there are subtle difference in
definitions between the subgroup version and the group homomorphism version.
A lot of representation features (e.g.\ dimension data) are closely related to the local
conjugacy. So the question at the beginning is closely related  to the
difference between the local and global conjugacy issue. $GL (n)$ case is simple since
by the character theory, local conjugacy implies global conjugacy. For general $G$,
this is not the case. When $H$ is finite, it is much easier to find the example of local conjugacy
without global conjugacy (\cite{Larsen94}, \cite{Larsen96}). When $H$ is connected, some pure
(but not so trivial) representation reasons
will tell. For more discussion, see Section ~\ref{S:2}.

\medskip

We are interested in such question for several reasons. One is related to
the multiplicity. When $G = \mGL (n)$,
local conjugacy implies global conjugacy, and then
this reflects the famous multiplicity one property for $\mGL (n)$.
For $G = \mPGL (n)$, local conjugate subgroups might not be conjugate, and this
reflects the failure of multiplicity one property
in $\mSL (n)$ (\cite{Blasius94}, \cite{Lapid99}).
When $G = \mSO (2 N)$, first connected instance of local conjugate but not globally conjugate
subgroups were found, and this induces an example the failure of multiplicity for $\mSO (2 N)$
with some assumptions on functoriality (\cite{Wang2007}), and such example
is totally different the ones
before and hence answered Langlands' question in some sense (\cite{Wang2007}, \cite{FLN2010}).
In fact, when $G = \mSO (10)$,
$H = \mSO (5)$, $\rho = \Lambda^{2}$: $H \to G$ the exterior square representation,
 and $\rho' = \tau \circ \rho$ where $\tau$ be an outer automorphism
of $\mSO (10)$ (i.e., a conjugation by an element of $\mO (10)$ with determinant $-1$),
$\rho$ and $\rho'$ are locally conjugate but not globally conjugate
in image. Assuming sufficient functoriality, when starting from a tempered modular
$4$-dim $l$-adic Galois representation of full $\mGSp (4)$ type which is easily obtained
from the Tate module of a principally polarized generic Abelian surface,
we get a $10$-dim orthogonal $l$-adic
Galois representation, which is associated to a (stable) cusp form of $\mSO (10)$
of multiplicity $> 1$. in fact, we found several families of such example. The later work
is to remove assumptions as much as we can.

% Main Results

According our philosophy, local conjugacy without global conjugacy in image will
leads to the failure of multiplicity one with sufficient functoriality condition.
For local conjugacy VS global conjugacy in group homomorphism form, there
are some subtle things, and it seems more strong assumption is needed to get
the failure of multiplicity one.

\medskip

This paper is to gather most of our local-global results in purely a representation way.
For discussions related to automorphic form, see \cite{Wang2011}. Throughout, unless specified,
a homomorphism means a group homomorphism
of complex algebraic groups, and a representation means a finite dimensional complex
representation.

\medskip

\begin{main} \label{TM:A}
\textnormal{(a)} Let $G = \mGL (N, \bC)$, $\mSL (N, \bC)$, $\mSp (2 N, \bC)$, $\mSO (2 N + 1, \bC)$
or $O (N, \bC)$, or
the group isogenous to above (except for $\mO (2 N, \bC)$), $H$ a connected
complex reductive group, $\rho, \rho'$: $H \to G$ be two homomorphisms.
If $\rho$ and $\rho'$ are locally conjugate, then they are globally conjugate.

\medskip

\textnormal{(b)} Let $G = \mSO (2 N, \bC)$ or its isogenous form,
and $H$, $\rho, \rho'$ as \textnormal{(a)}.
If $\rho$ and $\rho'$ are locally conjugate, then they either are globally conjugate,
or differ by an outer automorphism.

\end{main}

\bigskip

\emph{Remark}: Two connected algebraic groups $G$, $G'$ are said to be \emph{isogenous},
if there is a connected algebraic group $G''$ together with two finite central isogenies
$G'' \to G$ and $G'' \to G'$. We also say that $G'$ is \emph{an isogenous form} of $G$.
For explicit definitions, see Section ~\ref{SS:203}.

\medskip

This result for $G$ classical is well known in representation theory,
at least at the time of Dynkin
(\cite{Dynkin}, \cite{Dynkin2}). For completeness we include all proof in this paper.

\medskip

When $G = \mSO (2 N, \bC)$, or its isogenous form, if $\rho, \rho'$: $H \to G$ are locally
conjugate, then they must either be globally conjugate or
differ by an outer automorphism.
So in such case, we need just
study the case $\rho$ and $\rho' = \tau \circ \rho$ where $\tau$ is an outer automorphism
of $G$ which is induced by a conjugation of an element of $\mO (2 N)$ of determinant $-1$.
Moreover, explicit analysis on multiplicities enable us to just focus on the case $G = \mSO (2 N)$.

\medskip

\begin{main} \label{TM:B}
Let $G = \mSO (2 N, C)$, $H$ a connected reductive group,
$\rho$: $H \to G$ be a homomorphism,
and $\rho' = \tau \circ \rho$ where $\tau$ be an outer automorphism of $G$.
Then we have the following:

\medskip

\textnormal{(a)} $\rho$ and $\rho'$ are locally conjugate if and only if $\rho$ has weight $1$, namely,
the restriction of $\rho$ to the maximal torus has trivial representation as its subrepresentation.

\medskip

\textnormal{(b)} $\rho$ and $\rho'$ are globally conjugate if and only if $\rho$ has an odd dimensional
orthogonal subrepresentation.

\medskip

\textnormal{(c)} Assume moreover that $\rho$ is injective.
$\rho$ and $\rho'$ are globally conjugate in image if and only if some automorphism of
$H$ lifts to an outer automorphism of $\mSO (2 N)$, i.e., a conjugation by an element
of $\mO (2 N)$ of determinant $-1$.
\end{main}

\medskip

% LFMO-special: What a cute name!

We say that a homomorphism $\rho$: $H \to \mSO (2 N, \bC)$
is \emph{LFMO-special}, if, for
 $\rho' = \tau \circ \rho$ for an outer automorphism $\tau$ of $\mSO (2 N)$,
the following conditions for $\rho$ and $\rho'$ hold.

\medskip

(1) $\rho$ is \emph{essential}, namely, the image of $\rho$ is not contained in any parabolic
subgroup of $\mSO (2 N)$. Or equivalently, as a representation $\rho$
has no nontrivial totally isotropic subrepresentation.

\medskip

(2) $\rho$ and $\rho'$ are locally conjugate.

\medskip

(3) $\rho$ and $\rho'$ are not globally conjugate in image.

\medskip

In this case, we also say that the $2 N$-dimensional representation induced by $\rho$
is LFMO-special.

\medskip

If $\rho$ is LFMO-special, then $\rho$ and $\rho'$ will be ``leading to failure of
multiplicity one''. That is also why we call this name. For reasons, see Theorem 3.2 of
\cite{Wang2011}.

\medskip

\begin{main} \label{TM:C}
Let $H$ be a complex connected reductive group and $\rho: H \to G = \mSO (2N, \bC)$
a homomorphism. Then $\rho$ is LFMO-special if and only if:

\medskip

\textnormal{(1)} The representation space of $\rho$ decomposes as a direct sum
of inequivalent even dimensional orthogonal subrepresentations.

\medskip

\textnormal{(2)} $\rho$ has weight $1$.

\medskip

\textnormal{(3)} If an automorphism $\phi$ of $\rho (H)$
lifts to an automorphism of $G$, then it lifts to an inner automorphism
of $G$.

\medskip

Moreover, if $\rho'$ is quasi-equivalent to the representation induced by $\rho$,
then $\rho'$ is LFMO-special if and only if so is $\rho$.

\end{main}

\medskip

We say that $\rho$ is \emph{stable} if the $2 N$-representation induced
by $\rho$ (also called $\rho$) is irreducible, .
Let $H = T \times H_{1} \times \ldots \times H_{r}$, where $T$ is a complex torus,
$H_{i}$ are simple Lie groups and $\rho$ is a representation of $H$.
Then by basic representation theory,
we have $\rho = 1 \otimes \bigotimes_{i} \rho_{i}$
where $\rho_{i}$ is a representation of $H_{i}$, either
symplectic or orthogonal. For general connected $H'$,
$\rho'$ must be quasi-equivalent to some $\rho$ of $H$ of the above case.
Moreover, when $\rho$ is a stable LFMO-special representation,
$\rho'$ must have weight $1$, and hence $\rho'$ must kill
the center of $H'$. Thus $\rho'$ is quasi-equivalent to
$\rho$ of $H$ for some $H$ is semisimple of adjoint type.
So the study of stable LFMO-special representation can be reduced
to the case when $H$ is semisimple of adjoint type.

\medskip

\begin{main} \label{TM:D}
Let $H = T \times H_{1} \times \ldots \times H_{r}$,
$\rho = \bigotimes_{i} \rho_{i}$
be an even dimensional orthogonal representation of a complex connected
reductive group where $T$ is a complex torus, $H_{i}$
is a simple Lie group of adjoint type, and $\rho_{i}$
is self-dual and irreducible for each $i$.

\medskip

Then $\rho$ is LFMO-special if and only if $\rho (H)$ has trivial center,
and one of the following cases happens:

\medskip

\textnormal{Case (1)}: Exactly one $\rho_{i}$ is even dimensional.
In this case such $\rho_{i}$ is orthogonal and LFMO-special.

\medskip

\textnormal{Case (2)}: Exactly two, say $\rho_{i}$ and $\rho_{j}$,
are even dimensional.
In this case the following must be excluded:
$H_{i}$ and $H_{j}$ are isogenous, and $\rho_{i}$
and $\rho_{j}$ are quasi-equivalent, and moreover the dimension
of $\rho_{i}$ is $\equiv 2 \pmod{4}$.

\medskip

\textnormal{Case (3)}: At least three $\rho_{i}$s are even dimensional.
\end{main}

\medskip

\emph{Remark}: We say that two representations
$\rho$ and $\rho'$ of $H$ and $H'$
are \emph{quasi-equivalent} if and only if there are
finite central isogenies $\iota: H'' \to H$ and
$\iota: H'' \to H'$ such that as representations of $H''$,
$\rho \circ \iota$ and $\rho' \circ \iota'$
are conjugate.

\medskip

In particular, we have the following corollary,
which shows up also in \cite{Wang2007} and \cite{Wang2011}.

\medskip

\begin{mainc} \label{TM:E}
Let $H = H_{1} \times H_{2} \times \ldots \times H_{r}$ be a
semisimple Lie group of adjoint type over $\mathbb{C}$ with simple
factors $H_{j} (1 \leq j \leq r)$, among which at least one has even
rank. Let $V = \mathfrak{h}_{1} \otimes \mathfrak{h}_{2} \otimes
\ldots \otimes \mathfrak{h}_{r}$, $\kappa_{i}$ the Killing form of
$\mathfrak{h}_{i} = {\rm Lie} (H_{i})$, and $\kappa = \kappa_{1}
\otimes \kappa_{2} \otimes \ldots \otimes \kappa_{r}$.

$$
\rho = {\rm Ad}_{H_{1}} \otimes {\rm Ad}_{H_{2}} \otimes \ldots
\otimes {\rm Ad}_{H_{r}}: H \to G = \mSO (V, \kappa) \cong \mSO (2 N,
\mathbb{C})
$$

Then $\rho$ is LFMO-special if and only if one of the following happens:

\medskip

\textnormal{(a)} Only one of $H_{j}$ is of even rank, and
$\mathfrak{h}_{j}$ is $A_{4 n} (n \ge 1), B_{2 n} (n \ge 1), C_{2 n}
(n \ge 1), E_{6}, E_{8}, F_{4}, G_{2}$.

\medskip

\textnormal{(b)} Exactly two of them, say $H_{j}$ and $H_{k}$, are
of even rank, and either $H_{j} \not\cong H_{k}$, or $H_{j} \cong
H_{k}$ with $4 | {\rm dim} H_{j}$.

\medskip

\textnormal{(c)} At least three of them are of even rank.

\medskip

In this case $\rho$ will lead to the failure of multiplicity one for
$\mSO (2 N)$.
\end{mainc}

\emph{Remark}: For condition (b), $4 | {\rm dim} H_{j}$ if and only
if $\mathfrak{h}_{j}$ is of type $A_{2 n}, B_{4 n}, C_{4 n}$, $D_{4 n},
E_{8}, F_{4}$, as
\begin{align}
{\rm dim} A_{2 n} &= 4 n (n + 1), & {\rm dim} B_{2 n} &= 2 n (4 n +
1) & {\rm dim} C_{2 n} &= 2 n (4 n + 1), \notag \\
{\rm dim} D_{2 n} &= 2 n (4 n - 1), & {\rm dim} E_{6} &= 78 & {\rm
dim} E_{8} &= 248 \notag \\
{\rm dim} F_{4} &= 52 & {\rm dim} G_{2} &= 14 && \notag
\end{align}

% Organization of the article

The article is organized as following. In Section ~\ref{S:2},
we will focus on the structure theory of self-dual representations
of connected reductive groups, and prove Theorem ~\ref{TM:A}. In Section ~\ref{S:3},
we will study the $\mSO (2 N)$ case, and finally prove all other theorems.
Some parts are listed for the purpose of clarification
of notations and concepts, and experts
can skip them.

\medskip

% Acknowledge

% This note is completed while I was visiting the University of Science and Technology of China,
% and here I express my thanks to them to provide a good environment for work.
% Also in these year, we benefit from series of important discussions with D.~
% Ramakrishnan these years. Also, we are grateful for
% helpful comments from Binyong Sun during the writing of this
% article. Also, we also thank the Morningside Center of Mathematics in Beijing to provide a
% very good environment for the work.

\textbf{Acknowledgement}: In these year, we benefit from series of important discussions with D.~
Ramakrishnan these years. Also, we are grateful for
helpful comments from Binyong Sun during the writing of this
article. We also thank Prof. Yi Ouyang for inviting me to visit 
the University of Science and Technology of China as I finished this note there.

\bigskip

\section{Structure Theory} \label{S:2}

% To insert the summary of this section

\subsection{Preliminaries} \label{SS:201}

\ssmedskip

First let's have some notations and preliminary discussions. Let $G$ and $H$ be complex reductive groups.
Recall that two homomorphisms $\rho$ and $\rho'$ are \emph{locally conjugate} if for
$h \in H$ in a (Zariski) dense open subset in $H$, $\rho (h)$ and $\rho' (h)$ are conjugate in $G$;
$\rho$ and $\rho'$ are \emph{globally conjugate} if $\rho$ and $\rho'$ differ
by a conjugation by an element in $G$, i.e., there exists a $g \in G$ such that $\rho' (x) = g \rho (x) g^{-1}$
for all $x \in H$; And moreover, $\rho$ and $\rho'$ are \emph{globally conjugate in image} if $\rho (H)$
and $\rho' (H)$ are conjugate in $G$, i.e., there is a $g \in G$ such that $\rho' (H) = g \rho (H) g^{-1}$.

\medskip

Of course if $\rho$ and $\rho'$ are globally conjugate, then they are locally conjugate. For $G = \mGL (n)$,
the converse is also true. This is just the character theory. For other $G$, this is in general not the case.

\medskip
% Precise Language and examples

Now explain the concept of ``globally conjugate in image''. It arise from the following subtle
situation. Given $\rho$ and $\rho'$ as before, if $\rho (H)$ and $\rho' (H)$
are conjugate in $G$ (and hence also called globally conjugate),
then we don't necessarily have $\rho$ and $\rho'$ are globally conjugate.
In fact, we have the following easy lemma:

\medskip

\begin{lemma} \label{T:201}
Let $\rho$ and $\rho'$: $H \to G$ as before. If they are globally conjugate in image,
then there is a $\lambda \in {\rm Aut} (H)$ such that $\rho$ and $\rho' \circ \lambda$
are conjugate.
\end{lemma}

\medskip

\[
\begin{tikzcd}
H \arrow{r}{\rho} \arrow{d}{\lambda} &G \arrow{d}{c_{g}} \\
H \arrow{r}{\rho} &G
\end{tikzcd}
\]

\medskip

\emph{Proof}: By modifying by a conjugation, may
assume that $\rho (H)$ and $\rho' (H)$ have the same image. Then $\rho$ and $\rho'$
must differ by a conjugation.
\qedsymbol

\medskip

To measure the distance among these concepts, we introduce, like in \cite{Wang2011},
the multiplicity numbers.

\medskip

\emph{Definition--Multiplicities}: We denote $M (\rho; G)$ the set of equivalent classes of $\rho'$:
$H \to G$ which is locally conjugate to $i$, modulo the global conjugacy, and
$M' (\rho; G)$ the set of equivalent classes of $\rho'$:
$H \to G$ which is locally conjugate to $i$, modulo the global conjugacy in image.
Moreover, denote $m (\rho; G)$ and $m' (\rho; G)$ as the cardinalities of $M (\rho; G)$
and $M' (\rho; G)$ respectively.

\medskip

Here are examples and some of the known facts. Unless specified, $H$ denotes a complex reductive group.

\medskip

\emph{Example}: When $G = \mGL (n, \bC)$ and $H$ reductive and not necessarily connected,
 $m (\rho; G) = m' (\rho; G) = 1$.

\medskip

\emph{Example}: When $G = \mSL (N, \bC), \mSp (2 N, \bC), \mO (N, \bC), \mSO (2 N + 1, \bC)$,
and $H$ reductive, not necessarily connected, then $m (\rho, G) = m' (\rho, G) = 1$.
This is well known for a long time. To make this note complete, we include
the whole purely algebraic proof in Section ~\ref{SS:202}.

\medskip

\emph{Fact}: As $M' (\rho; G)$ is a quotient set of $M (\rho; G)$, $m' (\rho; G) \leq m (\rho; G)$
with the equality holds if ${\rm Aut} (H) = {\rm Inn} (H)$, i.e., $H$
has no outer automorphism.

\medskip

\emph{Example}: When $G = \mSO (2 N, \bC)$, then $m (\rho; G)$ and $m' (\rho; G)$ are
$1$ or $2$. One purpose of this note is to discuss this explicitly when $H$ is connected.
In particular, when $G = \mSO (10, \bC)$, $H = \mSO (5, \bC)$ and $\rho = \Lambda^{2}$, the
exterior square map, then $m' (\rho; G) = 2$. This is the simplest case when $H$ is connected and
$m' (\rho; G) = 2$.

% The reason is to be added.

\medskip

\emph{Example}: When $G = \mPGL (N, \bC)$, $m (\rho; G), m' (\rho; G)$ might not
be $1$ (\cite{Blasius94}, \cite{Lapid99}).
However, when $H$ is connected,
$m (\rho; G) = m' (\rho'; G) = 1$. See Section ~\ref{SS:203}.
.

\medskip

\emph{Examples}: When $G = F_{4}$, $E_{6}$, $E_{7}$ or $E_{8}$ (any isogenous form),
$H$ is finite, $m (\rho; G)$
and $m' (\rho; G)$ might be also greater than $1$ (see \cite{Larsen94}, \cite{Larsen96}).

\medskip

\emph{Example}: When $G = E_{8} (\bC)$ and $H = \mPGL (3, \bC)$,
there is an embedding $\rho$: $H \hookrightarrow G$
such that $m (\rho; G) = 2$ (\cite{Seitz87}, \cite{Seitz91}, and also Subsection ~\ref{SS:302})

\medskip

\emph{Fact}: for general $G$ and reductive $H$, we have $m (\rho; G)$ and $m' (\rho; G)$ are finite
(\cite{Dynkin}, \cite{Dynkin2}, \cite{Wang2011}, \cite{A-Y-Y2010}, \cite{Yu}).

\bigskip

\subsection{Structure Theory: Self-dual Representations} \label{SS:202}

\ssmedskip

In this part, we study the case when $G$ is classical and $H$ is reductive, not necessarily connected.
The main results of this part are Theorem ~\ref{T:202}, Corollary ~\ref{T:204},
Theorem ~\ref{T:205}. For review of the concepts,
see the discussion before the proof of Theorem ~\ref{T:202}.

% $O (2 N)$, $SO (2 N + 1)$ and $Sp (2 N)$ structure theorem

\medskip

\begin{theorem} \label{T:202}
Let $G = \mO (N, \bC)$, $\mSO (2 N + 1, \bC)$ or $\mSp (2 N, \bC)$, $H$ complex reductive, and
$\rho, \rho'$: $H \to G$ two homomorphisms. Then if $\rho, \rho'$
are locally conjugate, then they are globally conjugate. Hence,
$m (\rho; G) = m' (\rho; G) = 1$. In particular, two equivalent orthogonal / symplectic
representations of $H$ are isometric.
\end{theorem}

\medskip

Before we start the proof, we list something here.

\medskip

\begin{proposition} \label{T:203}
When $G = \mGL (N, \bC), \mSL (N, \bC)$, $H$ reductive and $\rho: H \to G$ a rational representation,
then $m (\rho; G) = m' (\rho; G) = 1$.
\end{proposition}

\medskip

\qedsymbol

\medskip

\begin{corollary}  \label{T:204}
When $G = \mSO (2 N, \bC)$, $H$, $\rho$ as before, then $m (\rho; G)$ and $m' (\rho; G)$ are $1$ or $2$.
In particular, if $\rho$ and $\rho'$ are locally conjugate, then either they are globally conjugate,
or they differ by an outer automorphism of $\mSO (2 N, \bC)$.
\end{corollary}

\medskip

\emph{Proof by using Theorem ~\ref{T:202}}:

\medskip

Let $\tilde{G} = \mO (2 N, \bC)$, $I: G \hookrightarrow \tilde{G}$ the natural embedding,
and $\rho': H \to G$ be any homomorphism which is locally conjugate
to $\rho$. Then $I \circ \rho$ and $I \circ \rho'$ are also locally conjugate.
By Theorem ~\ref{T:202}, $I \circ \rho$ and $I \circ \rho'$ are globally conjugate
in $\tilde{G}$. Let $g \in \mO (2 N, \bC)$ be such that $I \rho' (x) = g (I \circ \rho (x)) g^{-1}$
for all $x \in G$. Then we have ${\rm det} (g) = \pm 1$. Let $c_{g}$ be the conjugation by $g$.
Then $c_{g}$ stabilizes $I (G)$ and hence leads to an automorphism of $G$, also called $c_{g}$,
and $\rho' = c_{g} \circ \rho$. If ${\rm det} (g) = 1$, then $c_{g}$ is inner;
If ${\rm det} (g) = -1$, then $c_{g}$ is an outer automorphism. Hence we have the corollary.

\qedsymbol

\medskip

Now we analyze the structure of the self-dual representation. First recall some
notations and concepts.

\medskip

Recall that given a vector space $W$ together with a bilinear form $B$
(not necessarily with a group action), a subspace $W'$
is said to be \emph{non-degenerate} (resp.\ \emph{totally isotropic})
if when $W'$, together with the bilinear form $B |_{W'}$, is non-degenerate
(resp.\ totally isotropic). (Recall that a space $W$ with a bilinear form $B$
is said to be \emph{totally isotropic} if $B (v, w) = 0$ for all $v, w \in W$.)

\medskip

Recall that,
given a finite-dimensional complex representation $(\rho, V)$ of a group $H$,
the \emph{contragredient} $(\rho^{\vee}, V^{\vee} = V^{*})$ of $\rho$ is
defined as, $V^{*} = {\rm Hom}_{\bC} (V, \bC)$ and $G$ acts as the following:
$(\rho^{\vee} (g) (\lambda)) (v) = \lambda (\rho (g^{-1}) v)$
for $\lambda \in V^{*}, v \in V, g \in H$.
A representation $(\rho, V)$ of $H$ is said to be \emph{self-dual} if $\rho^{\vee}$
is equivalent to $\rho$, which is equivalent to say that, there is a non-degenerate
bilinear form $B: V \times V \to \bC$ which is $H$-equivariant, i.e.,
$B (\rho(g) v, \rho (g) w) = B (v, w)$ for $v, w \in V$.
Moreover, it is well known that, if $\rho$ is irreducible and self-dual,
then such $B$ must be
unique up to scalar (Schur's Lemma, see Lemma ~\ref{T:206}),
and moreover must be symmetric or alternating.
We also say that in this case $\rho$ is \emph{of orthogonal type} (resp.\ \emph{of symplectic type})
or just \emph{orthogonal} (resp.\ \emph{symplectic})
if $B$ is symmetric (resp.\ alternating). A \emph{$H$-invariant subspace}
or \emph{$H$-subspace} of $V$ is a
subspace which is stable under the action of $H$.

\medskip

Note that when we talk about the \emph{orthogonal/symplectic structure} of a self-dual
representation $(\rho, V)$ of $H$ we also means together a non-degenerate $H$-invariant
symmetric/alternating bilinear form $B$, which is not unique in general but is unique
up to scalar when $\rho$ is irreducible. In this case, $\rho$ is induced
by a homomorphism $\rho_{0}$: $H \to G$ where $G = \mO (N), \mSp (2 N)$ (or $\mSO (N)$
if the image of $\rho$ happens to lies inside $\mSL (V)$) through
the standard representation $G \overset{\cong}{\to} \mO (V, B) \hookrightarrow \mGL (V)$
where $\mO (V, B) := \{ g \in \mGL (V) \mid B (\rho (g) v, \rho (g) w) = B (v, w) \,\forall v, w \in V \}$.
$\rho_{0}$ or $\rho$ is called \emph{essential} if any totally isotropic
$H$-invariant subspace of $V$ must be trivial, and $\rho_{0}$ or $\rho$
is called \emph{stable} if $\rho$ is irreducible. Of course if $\rho$ or $\rho_{0}$
is stable then it is essential. Moreover, it is well known that $\rho_{0}$ is essential if and only if
the image of $\rho_{0}$ is not contained in any proper parabolic subgroup of $G$.
So the concept ``essential'' should be the analogue of ``irreducible'' to the case
$G$ classical.

\medskip

\emph{Proof of Theorem ~\ref{T:202}}:

\medskip

Let $V$ and $V'$ be finite dimensional self-dual complex representations of $H$,
and $B$ and $B'$ are orthogonal/symplectic structures on $V$ and $V'$ induced by $\rho$ and $\rho'$
respectively. As $\rho$ and $\rho'$ are locally conjugate, $V$ and $V'$ must be equivalent as
$H$-representations (Proposition ~\ref{T:203}). Now our theorem follows from the following claim:
There exists an isomorphism $T: V \to V'$ which is $H$-equivariant and $B = {B'}^{T}$
where ${B'}^{T} (u, v) = B' (T u, T v) \,\forall u, v \in V$, i.e, $T$ is an $H$-equivariant isometry.

\medskip

Now we prove the claim:
Let $T_{0}: V \to V'$ be an arbitrary $H$-equivariant isomorphism. It is easy
to see that ${B'}^{T_{0}}$ is also a non-degenerate $H$-invariant bilinear form
on $V$. Through $B$ and $B'^{T_{0}}$ we get two $H$-representation isomorphisms
from $V$ to $V^{*}$. Thus there is a $L \in \mGL (V)$ which is $H$-equivariant
such that the following diagram commutes.

\medskip

\[
\begin{tikzcd}
V \arrow{r}{i_{B'^{T_{0}}}} \arrow{d}{L} & V^{*} \arrow{d}{=} \\
V \arrow{r}{i_{B}} & V^{*}
\end{tikzcd}
\]
where $i_{B}$: $V \to V^{*}$ is the isomorphism induced by $B$, i.e.,
$i_{B} (w) (v) = B (v, w)$ $\, \forall v, w \in V$.
In particular, ${B'}^{T_{0}} = B^{L}$. Hence $B = B'^{T_{0} \circ L^{-1}}$,
and $T = T_{0} \circ L^{-1}$ is what we want.

\qedsymbol

\medskip

\begin{theorem} \label{T:205}
Let $V$ be a finite dimensional self-dual complex
representation of a complex reductive group $H$ with $B$ a non-degenerate
$H$-invariant bilinear form on $V$ which is either symmetric or alternating.
Then we have the following.

\medskip

\textnormal{(1)} $V$ is a semisimple representation of $H$, i.e.,
$V$ is a direct sum of irreducible representations.

\medskip

\textnormal{(2)} Each irreducible $H$-subspace of $V$ is
either non-degenerate or totally isotropic.

\medskip

\textnormal{(3)} Let $V_{\sigma}$ be the \emph{$\sigma$-isotypical}
component of $V$, i.e., the sum of all irreducible subspace
equivalent to $\sigma$, and let $W_{\sigma} = V_{\sigma}$ if
$\sigma$ is self-dual, and $V_{\sigma} + V_{\sigma^{\vee}}$ if
$\sigma$ is not self-dual. Then $V$ is the direct sum of
$W_{\sigma}$, and the direct sum is orthogonal, and each
$W_{\sigma}$ is non-degenerate.

\medskip

\textnormal{(4)} If $\sigma$ is not self-dual, then $W_{\sigma} =
V_{\sigma} + V_{\sigma^{\vee}}$ gives rises to a complete
polarization of $W_{\sigma}$. If $\sigma$ is self-dual, then
$W_{\sigma} = V_{\sigma} \cong W_{\sigma, 0} \otimes W'_{\sigma}$
where $W_{\sigma, 0} \cong \sigma$ and $H$ acts on $W'_{\sigma}$
trivially. Moreover, we can endow $W'$ an orthogonal/symplectic
structure such that the isomorphism is also an isometry.

\medskip

\end{theorem}

\medskip

\begin{lemma} \label{T:206}
\textnormal{\textbf{(Schur's Lemma)}}
Let $V, W$ be two complex finite-dimensional irreducible representations
of a group $H$, and $T \in {\rm Hom}_{H} (V, W)$. Then

\medskip

\textnormal{(1)} $T$ is either $0$ or an isomorphism.

\medskip

\textnormal{(2)} If $V = W$, then $T = c \,{\rm id}_{V}$, a scalar multiplication.

\medskip

\textnormal{(3)} If $B, B'$ are two non-degenerate $H$-invariant bilinear forms on $V$,
and $B$ is non-degenerate, then $B' = c B$ for some $c \in \bC$.
\end{lemma}

\medskip

\qedsymbol

\medskip

\emph{Proof of Theorem ~\ref{T:205}}:

(1) follows from Proposition ~\ref{T:203}.

\medskip

(2): for each irreducible $H$-subspace $W$, $W \cap W^{\bot}$ is
also $H$-stable, where $W^{\bot} = \{ w' \in V \mid B (w, w') = 0 \quad \forall w \in W\}$.
As $W$ is irreducible, $W \cap W^{\bot} = W$ or
$0$. Hence the assertion.

\medskip

(3): First, let $W$ and $W'$ be two irreducible subspaces of $V$,
and assume that $W'$ is not congruent to $W^{\vee}$, then we claim
that $W$ and $W'$ are orthogonal. In fact, as $W'$ is not equivalent
to $W^{*} \cong W^{\vee}$, the map $W' \to W^{*}$, $w \mapsto
B(\,\*\,,  w)$, which is $H$-equivariant, must be $0$ by
Lemma ~\ref{T:206}. Hence the claim.

\medskip

Thus, $V_{\sigma}$ and $V_{\sigma'}$ are orthogonal if $\sigma'$ is
not equivalent to $\sigma^{\vee}$, and $W_{\sigma}$ and
$W_{\sigma'}$ are orthogonal if $\sigma'$ is not equivalent to
$\sigma$ or $\sigma^{\vee}$.

\medskip

Thus $V$ is the orthogonal direct sum of $W_{\sigma}$ when $\sigma$
runs through all pairs $\{\sigma, \sigma^{\vee}\}$ of irreducible
representations of $H$. Moreover, all $W_{\sigma}$ are
non-degenerate since $V$ itself is non-degenerate.

\medskip

(4): For $\sigma$ not self-dual, $V_{\sigma}$ and
$V_{\sigma^{\vee}}$ form a complete polarization of $W_{\sigma}$
since each irreducible subspace of type $\sigma$ must be totally
isotropic.

\medskip

For $\sigma$ self-dual, want to analyze $W_{\sigma} = V_{\sigma}$.
As it is also non-degenerate from (3), may assume $V = V_{\sigma}$,
with the non-degenerate $H$-invariant bilinear form $B$, either
symmetric or alternating.

\medskip

Let $W = W_{\sigma, 0}$ be an irreducible representation of $H$ of
type $\sigma$, with a non-degenerate $H$-invariant
symmetric/alternating form $B_{0}$. Thus $V \cong W \otimes W'$
as $H$-representations as $V = V_{\sigma}$ is $\sigma$-isotypical,
where the action of $H$ on the right side is on the first factor. We
identify both sides for such isomorphism.

\medskip

We claim that we can endow a non-degenerate bilinear form $B'$
on $W'$, such that
\[
B (u \otimes u', v \otimes v') = B_{0} (u, v) B' (u',
v')
\]
for each $u, v \in W$ and $u', v' \in W'$.

\medskip

Granting this claim, then set $W'_{\sigma} = W'$, and moreover
$B'$ is symmetric or alternating according to the types of
$B$ on $V \cong W \otimes W'$ and $B_{0}$ on $W$. Hence
the second assertion of (4) follows. The rest assertions follow
also.

\medskip

Now come to the claim. First, for each pair $(u, v) \in W^{2}$ that
$B_{0} (u, v) \ne 0$, define a bilinear form $B'_{(u, v)}$
on $W'$ such that
\[
B (u \otimes u', v \otimes v') = B_{0} (u, v) B'_{(u,
v)} (u', v')
\]
for each $u', v' \in W'$. The subscript for $B'_{(u, v)}$ is present since
the definition depends on the choice of $(u, v)$, and we will
see that finally the subscript can be dropped.

\medskip

We want to prove: (a) $B'_{(u, v)}$ is independent of the
choice of $(u, v)$. (b) $B' = B'_{(u, v)}$ is
non-degenerate, and the claim holds for such $B'$.

\medskip

Assume that $B'_{(u_{0}, v_{0})} (u', v') \ne 0, \,u_{0}, v_{0}
\in W,\, u', v' \in W'$. Such $u_{0}, v_{0}, u', v'$ exist as
$B$ is non-degenerate and there are $u, v \in W$, $u', v' \in
W'$ such that $B (u \otimes u', v \otimes v') \ne 0$.

\medskip

Fix $u'$ and $v'$, and put $\tilde{B}_{0}$ a bilinear form on $W$ such that $B
(u \otimes u', v \otimes v') = \tilde{B}_{0} (u, v)
B'_{(u_{0}, v_{0})} (u', v')$. Then $\tilde{B}_{0}$ is
also $H$-invariant. As $W$ is irreducible, then $\tilde{B}_{0}
= C B_{0}$ for some complex number $C$ by Lemma ~\ref{T:206}. As
$\tilde{B}_{0} (u_{0}, v_{0}) = B_{0} (u_{0}, v_{0})$,
$\tilde{B}_{0} = B_{0}$. Thus we have $B (u \otimes
u', v \otimes v') = B_{0} (u, v) B'_{(u_{0}, v_{0})} (u',
v')$ for ALL $u, v \in W$. So for all $u, v \in W$ with $B_{0}
(u, v) \ne 0$, $B'_{(u, v)} (u', v') = B'_{(u_{0}, v_{0})}
(u', v') \ne 0$. Since $B'_{(u, v)} (u', v')$, as a quotient $B (u \otimes u', v \otimes v')$
by $B_{0} (u, v)$,
is a rational function
of $W^{2} \times W'^{2}$, $B'_{(u, v)} (u, v)$ must be independent of the
choice of $(u, v)$ and hence (a).

\medskip

Let $B' = B'_{(u, v)}$ for any choice of $(u, v)$.
Then $B'$ is a bilinear form
on $W'$, and then (b) and the claim follow easily.

\medskip

\qedsymbol

\medskip

For later use, we quote the following.

\begin{proposition} \label{T:207}
Let $(\rho, V)$ be a self-dual finite dimensional representation of a complex reductive
group $H$ with a $H$-invariant symmetric/alternating bilinear form $B$.
Then $\rho$ is essential if and only if $\rho$ is a direct sum of
inequivalent self-dual representations.
Let $V_{\sigma}$ be as in Theorem ~\ref{T:205} for each irreducible representation $\sigma$.
Then in this case,
either $V_{\sigma} = 0$ or $V_{\sigma}$ is irreducible and non-degenerate.
\end{proposition}

\medskip

\emph{Proof}: All notations are the same as in Theorem ~\ref{T:205}.

\medskip

(1) Assume that $\rho$ is essential. By Theorem ~\ref{T:205} (2), all irreducible
constituents $\sigma$ of $V$ must be self-dual, and by Theorem ~\ref{T:205}
(4), we have $V = \bigoplus_{\sigma} V_{\sigma}$,
while $V_{\sigma} \cong W_{\sigma, 0} \otimes W'_{\sigma}$ being the isometry,
where $\sigma$
runs through all irreducible self-dual constituents $\sigma$ of $V$. If for some $\sigma$,
$W'_{\sigma}$ has dimension $\geq 2$, then $W'$ must have an isotropic vector.
Then $W_{\sigma, 0} \otimes \bC w'$ is a totally
isotropic $H$-subspace of $V_{\sigma} = W_{\sigma, 0} \otimes W'_{\sigma}$.
Through the isometry, this gives rise to a nontrivial totally isotropic
$H$-subspace. This contradicts the assumption that $\rho$ is essential.
Hence, Each $W'_{\sigma}$ must be $1$-dimensional and $V_{\sigma}$ must be irreducible.
In particular, $V$ is a direct sum of inequivalent irreducible non-degenerate $H$-subspaces
(Theorem ~\ref{T:205} (3)).

\medskip

(2) Assume now that $V$ is a direct sum of inequivalent irreducible
non-degenerate $H$-subspaces, and $W$ a totally isotropic $H$-subspace.
Want to prove that $W = 0$. Otherwise, let $W'$ is an irreducible $H$-subspace of $W$,
and want to get the contradiction. By Theorem ~\ref{T:205},
$W' \subset W_{\sigma}$ for some self-dual $\sigma$. By the assumption,
for each constituent $\sigma$, $\sigma$ is self-dual and $W_{\sigma} = V_{\sigma}$
is itself irreducible. Thus $W' = V_{\sigma}$ is non-degenerate (Theorem ~\ref{T:205}),
and hence $W \supset W'$ can't be totally isotropic. Contradiction. Done.

\medskip

\qedsymbol

\bigskip

\subsection{About Isogenous Forms} \label{SS:203}

\ssmedskip

From now on we focus on the case when $H$ a connected reductive group, still,
$\rho$ and $\rho'$ are homomorphisms from $H$ to $G$. To prove Theorem ~\ref{TM:A},
we need to study the relations between the multiplicities
and the isogenous forms.
The main results are Theorem ~\ref{T:208} and Theorem ~\ref{T:214},
and at the end of this subsection we will prove
Theorem ~\ref{TM:A}.

\medskip

First some definitions. Unless specified, all groups involved in this section are connected complex reductive.
We say that $\pi: \tilde{G} \to G$ is a \emph{finite central isogeny}
if $\pi$ is a finite homomorphism
with its kernel central finite.
In this case, we also say that $\tilde{G}$ is a \emph{finite central isogenous cover}
(or just \emph{isogenous cover}) of $G$.
It is not hard to see that the composition of two finite central isogenies
of complex reductive groups is still a finite central isogeny.
Moreover, two connected groups $G$ and $G'$ are said to be \emph{isogenous}
or \emph{of the same isogenous form}
if they share a common finite central isogenous cover, i.e., there is a $G''$ such that
both $G'' \to G$ and $G'' \to G'$ are finite central isogenies. In this case we also say that $G'$
is an \emph{isogenous form} of $G$. It is known that two connected reductive groups are
isogenous if and only if they share the same Lie algebra.

\medskip

We say that $\tilde{\rho}: \tilde{H} \to \tilde{G}$
is a \emph{finite central isogenous cover} of $\rho: H \to G$ if the following diagram commutes
\[
\begin{tikzcd}
\tilde{H} \arrow{r}{\tilde{\rho}} \arrow{d}{\pi_{H}} & \tilde{G} \arrow{d}{\pi_{G}} \\
H \arrow{r}{\rho} & G
\end{tikzcd}
\]
where $\pi_{H}$ and $\pi_{G}$ are also finite central isogenies.
We say that $\rho: H \to G$ and $\rho_{1}: H_{1} \to G_{1}$
are \emph{isogenous} if they possess the same finite central isogenous covering.
In this case, $H$ and $H_{1}$, $G$ and $G_{1}$ are of the same isogenous form.

\medskip

\begin{theorem} \label{T:208}
Let $\rho: H \to G$ and $\rho': H' \to G'$ be two isogenous homomorphisms
of connected complex reductive groups.
Then $m (\rho; G) = m (\rho'; G')$ and $m' (\rho; G) = m' (\rho'; G')$.
\end{theorem}

\medskip

This theorem enable us to reduce the case to a convenient isogenous
form of $G$ in our discussion.
For the definition of multiplicities, see Section ~\ref{SS:201}.

\medskip

According to the definition, it suffices for us to prove the case when one of $\rho'$
and $\rho$ is a finite central isogenous cover of another.

\medskip

\begin{proposition} \label{T:209}
Let $H$ be a connected reductive group and $\pi: \tilde{G} \to G$ a finite central isogeny
of connected complex reductive algebraic groups. Given a homomorphism $\rho$, we can lift
$\rho$ to $\tilde{\rho}: H \to \tilde{G}$ such that $\rho = \pi \circ \tilde{\rho}$
if and only if we can lift $\rho|_{T}$ where $T$ is a maximal torus of $H$.

\medskip

\[
\begin{tikzpicture}{>=angle 90}
\matrix(a)[matrix of math nodes, row sep=3em, column sep=2.5em,
text height=1.5ex, text depth=0.25ex]
{& & \tilde{G} \\ T & H & G \\};
\path[right hook->, font=\scriptsize] (a-2-1) edge (a-2-2);
\path[dotted, ->, font=\scriptsize] (a-2-1.north east) edge node[above]{$\tilde{\rho}_{T}$} (a-1-3);
\path[->, font=\scriptsize] (a-1-3) edge node[right]{$\pi$} (a-2-3);
\path[dotted, ->, font=\scriptsize] (a-2-2.north east) edge node[above]{$\tilde{\rho}$} (a-1-3);
\path[->, font=\scriptsize] (a-2-2) edge node[below]{$\rho$} (a-2-3);
\end{tikzpicture}
\]

\medskip

In this case, such $\tilde{\rho}$ is unique.
\end{proposition}

\medskip

\emph{Remark}: In this lemma, we call such $\rho$ \emph{liftable}.

\medskip

\emph{Proof}:

\medskip

The only if part of the first assertion is not a problem.
Once the first assertion done, the uniqueness is also easy to see due to the connectedness of
$H$. Now we prove the if part. We proceed by cases.

\medskip

Case (1): $H$ is semisimple and simply connected. In this case,
such $\tilde{\rho}$ definitely exists, even without the assumption
about $\rho |_{T}$.

\medskip

Case (2): $H$ is semisimple. Let $\tilde{H}$ be the simply connected complex Lie group
isogenous to $H$, and $\pi_{H}: \tilde{H} \to H$.
\[
\begin{tikzpicture}{>=angle 90}[descr/.style={fill=white}]
\matrix(a)[matrix of math nodes, row sep=3em, column sep=2.5em,
text height=1.5ex, text depth=0.25ex]
{\tilde{T} & \tilde{H} & \tilde{G} \\ T & H & G \\};
\path[right hook->, font=\scriptsize]
(a-1-1) edge (a-1-2) (a-2-1) edge (a-2-2);
\path[dotted, ->, font=\scriptsize]
(a-1-2) edge node[above]{$\tilde{\tilde{\rho}}$} (a-1-3)
(a-2-2) edge node[above]{$\tilde{\rho}$} (a-1-3);
\path[->, font=\scriptsize]
(a-1-1) edge node[left]{$\pi_{T}$} (a-2-1)
(a-1-2) edge node[left]{$\pi_{H}$} (a-2-2)
(a-1-3) edge node[right]{$\pi$} (a-2-3)
(a-2-2) edge node[below]{$\rho$} (a-2-3);
\end{tikzpicture}
\]
Here $\tilde{T}$ is the inverse image of $T$ via $\pi_{H}$ which is again
a maximal torus of $\tilde{H}$.
From case (1), $\tilde{\tilde{\rho}}$ exists.

\medskip

To prove that $\tilde{\rho}$ exists it suffices to show that
 $\tilde{\tilde{\rho}} \left( {\rm Ker} (\pi_{H}) \right) = 1$. This follows as  since $\rho |_{T}$
is liftable, and hence $\tilde{\tilde{\rho}} \left( {\rm Ker} (\pi_{T}) \right)  = 1$, and moreover,
as $\pi_{H}$ is a finite central isogeny, its kernel must be semisimple and contained in
all maximal tori of $\tilde{H}$, and hence ${\rm Ker} (\pi_{T}) = {\rm Ker} (\pi_{H})$.

\medskip

Case (3): $H$ itself is a torus. This is trivial as $H = T$ and $\rho = \rho |_{T}$.

\medskip

As $H$ is connected, we deduce that $\tilde{\rho}$ is unique, for at least the last
three cases, and we will use this in the last case.

\medskip

Case (4): General case. $H$ connected reductive.

\medskip

\[
\begin{tikzpicture}{>=angle 90}[descr/.style={fill=white}]
\matrix(a)[matrix of math nodes, row sep=3em, column sep=2.5em,
text height=1.5ex, text depth=0.25ex]
{& \tilde{G} & \\ H_{0} & H & T_{0} \\};
\path[dotted, ->, font=\scriptsize]
(a-2-1) edge node[left]{$\tilde{\rho}_{H_{0}}$} (a-1-2)
(a-2-3) edge node[right]{$\tilde{\rho}_{T_{0}}$}(a-1-2);
\path[->, font=\scriptsize]
(a-1-2) edge node[left]{$\pi$} (a-2-2)
(a-2-1) edge node[below]{$\pi_{H_{0}}$} (a-2-2)
(a-2-3) edge node[below]{$\pi_{T_{0}}$} (a-2-2);
\end{tikzpicture}
\]

Let $H_{0} = (H, H)$ be the derived group of $H$, $T_{0} \subset T$ the maximal central
torus of $H_{0}$, and $\rho_{H_{0}} = \rho |_{H_{0}}$ and $\rho_{T_{0}} = \rho |_{T_{0}}$,
and let $T_{1} = H_{0} \cap T$ be also a maximal torus of $H_{0}$.
As $\rho |_{T}$ is liftable, so is $\rho |_{T_{1}}$ and $\rho_{T_{0}}$.
Applying Case (2) (while $H_{0} = (H, H)$ is semisimple) and (3), we get, both $\rho_{H_{0}}$ and
$\rho_{T_{0}}$ are liftable. Now define $\tilde{\rho}: H \to G$ as
\[
\tilde{\rho} (h_{0} t_{0}) = \tilde{\rho}_{H_{0}} (h_{0}) \tilde{\rho}_{T_{0}} (t_{0}) \qquad \forall h_{0} \in H_{0},
t_{0} \in T_{0}
\]
We claim that it is well defined, and is a homomorphism.

\medskip

If $h_{0} t_{0} = 1$, then $h_{0}, t_{0} \in Z (H_{0}) \subset T_{1}$, and in particular, lies in
$T$. Thus by the uniqueness of $\tilde{\rho}_{T}$, $\tilde{\rho}_{H_{0}}$
and $\tilde{\rho}_{T_{0}}$ and $\tilde{\rho}_{T_{1}}$ where $\tilde{\rho}_{T_{1}} = \tilde{\rho}_{T} |_{T_{1}}$
is the
lift of $\rho_{T_{1}}$, we have
\begin{align}
\tilde{\rho} (h_{0} t_{0}) &= \tilde{\rho}_{H_{0}} (h_{0}) \tilde{\rho}_{T_{0}} (t_{0})
& &\text{(Definition.)} \notag\\
&= \tilde{\rho}_{T_{1}} (h_{0}) \tilde{\rho}_{T_{1}} (t_{0})
& & \text{($h_{0} \in T_{1}$ and $T_{1} \subset H_{0}$)} \notag\\
&= \tilde{\rho}_{T_{1}} (h_{0} t_{0}) = 1 &&\text{$h_{0} t_{0} = 1$} \notag
\end{align}
As $H_{0}$ and $T_{0}$ commute, $\tilde{\rho}$ is well defined. Now it is routine to check that
$\tilde{\rho}$ is a homomorphism and is a lift of $\rho$, and finally is unique.

\qedsymbol

\medskip

\begin{lemma} \label{T:210}
Let $T$ be a complex torus.
Then $T$ has a \emph{generic point}, i.e., a point that
generates a (Zariski or topologically) dense subset of $T$.
the set of generic points of $T$ is (Zariski or topologically) dense.
\end{lemma}

\medskip

\qedsymbol

\medskip

\begin{lemma} \label{T:211}
Let $\rho, \rho': H \to G$ be two homomorphisms
of complex connected reductive groups.
Let $T$ be maximal torus of $H$. Then $\rho$ and $\rho'$ are locally conjugate,
if and only if $\rho |_{T}, \rho' |_{T}: T \to G$ are globally conjugate.
\end{lemma}

\medskip

\emph{Proof}:

\medskip

(a) If part: Assume that $\rho |_{T}$ and $\rho' |_{T}$ are globally conjugate.
From linear algebraic group theory, $H^{ss}$,
the set of semisimple points of $H$, is Zariski dense in $H$, and contains an open subset
of $H$, and moreover each semisimple element of $H$
is conjugate to en element of $T$. Thus for each $h \in H^{ss}$,
$\rho (h)$ and $\rho' (h)$ are conjugate in $G$, and
hence $\rho$ and $\rho'$ are locally conjugate.

\medskip

(b) Only if part: Assume that $\rho$ and $\rho'$ are locally conjugate.
By the definition, $\rho (h)$ and $\rho' (h)$ are conjugate in $G$
for $h$ in a (Zariski) dense open subset $W$ in $H$. Since $H^{ss}$
is Zariski dense in $H$, $W \cap G^{ss} \ne \emptyset$
and hence there is a maximal torus $T_{0}$ which intersects $W$. As $W$ is Zariski open
dense in $H$, $W \cap T_{0} \ne \emptyset$ must be also Zariski open dense in $T_{0}$.
From Lemma ~\ref{T:210}, there is a $t_{0} \in T_{0} \cap W$ which is generic.
Hence $\rho |_{T_{0}}$ and $\rho' |_{T_{0}}$ must be globally conjugate. The maximal tori $T$
and $T_{0}$ are conjugate (Well known!), and hence
$\rho |_{T}$ and $\rho' |_{T}$ must be globally conjugate.

\medskip

\qedsymbol

\medskip

\begin{corollary}  \label{T:212}
Let $H$, $G$, $\tilde{H}$, $\tilde{G}$ be connected complex reductive groups,
$\pi_{H}: \tilde{H} \to H$, $\pi_{G}: \tilde{G} \to G$ finite central isogenies,
$\rho, \rho': H \to G$ be homomorphisms. Assume

\medskip

\textnormal{(1)} $\rho$ and $\rho'$ are locally conjugate.

\medskip

\textnormal{(2)} $\rho$ lifts to $\tilde{\rho}$ (see the diagram).

\medskip

Then $\rho'$ also lifts to $\tilde{\rho}'$, and the lift is unique.

\[
\begin{tikzcd}
\tilde{H} \arrow{r}{\rho, (\text{?} \tilde{\rho}')} \arrow{d}{\pi_{H}}
&\tilde{G} \arrow{d}{\pi_{G}} \\
H \arrow{r}{\rho, \rho'} &G
\end{tikzcd}
\]
\end{corollary}

\medskip

\emph{Proof}: Let $\tilde{T}$ be a maximal torus of $\tilde{H}$
and $T = \pi_{H} (\tilde{T})$.

\medskip

Consider $\rho \circ \pi_{H}$. As $\rho$ lifts $\tilde{\rho}$, $\rho \circ \pi_{H}
= \pi_{G} \circ \tilde{\rho}$. Hence $\rho_{\tilde{T}} \overset{\Delta}{=}
(\rho \circ \pi_{H}) |_{\tilde{T}}$ is also liftable.

\medskip

As $\rho$ and $\rho'$
are locally conjugate, so is $\rho \circ \pi_{H}$ and $\rho' \circ \pi_{H}$.
Pick In particular,
$\rho_{\tilde{T}}$ and $\rho'_{\tilde{T}} \overset{\Delta}{=} \rho' \circ \pi_{H}$
are globally conjugate (Lemma ~\ref{T:211}). Thus,
$\rho'_{\tilde{T}}$ is also liftable.

\medskip

\[
\begin{tikzpicture}{>=angle 90}[descr/.style={fill=white}]
\matrix(a)[matrix of math nodes, row sep=3em, column sep=2.5em,
text height=1.5ex, text depth=0.25ex]
{\tilde{T} & \tilde{H} & \tilde{G} \\ & H & G \\};
\path[dotted, ->, font=\scriptsize]
(a-1-2) edge node[above]{$\tilde{\rho}'$} (a-1-3);
\path[right hook->, font=\scriptsize]
(a-1-1) edge (a-1-2);
\path[->, font=\scriptsize]
(a-1-1) edge (a-2-2)
(a-1-2) edge node[right]{$\pi_{H}$} (a-2-2)
(a-1-3) edge node[right]{$\pi_{G}$} (a-2-3)
(a-2-2) edge node[below]{$\rho'$} (a-2-3);
\end{tikzpicture}
\]

\medskip

Then by Proposition ~\ref{T:209}, $\rho' \circ \pi_{H}$ is liftable.
This gives rise to the existence of $\tilde{\rho}'$. The uniqueness now follows
easily from Proposition ~\ref{T:209}.

\medskip

\qedsymbol

\medskip

\begin{theorem}  \label{T:213}
Let $H$, $G$, $\tilde{H}$, $\tilde{G}$ be connected complex reductive groups,
$\pi_{H}: \tilde{H} \to H$, $\pi_{G}: \tilde{G} \to G$ finite central isogenies,
$\rho, \rho': H \to G$ be homomorphisms which lift to
$\tilde{\rho}, \tilde{\rho}'$ respectively.

\[
\begin{tikzcd}
\tilde{H} \arrow{r}{\rho, (\text{?} \tilde{\rho}')} \arrow{d}{\pi_{H}}
&\tilde{G} \arrow{d}{\pi_{G}} \\
H \arrow{r}{\rho, \rho'} &G
\end{tikzcd}
\]

\medskip

Then $\tilde{\rho}$ and $\tilde{\rho}'$ are locally conjugate / globally conjugate /
globally conjugate in image if and only if so are $\rho$ and $\rho'$.
\end{theorem}

\medskip

\emph{Proof}: We proceed in steps.

\medskip

Step 1: $\rho, \rho'$ are globally conjugate $\quad\Longrightarrow\quad$ $\tilde{\rho},
\tilde{\rho}'$ are globally conjugate.

\medskip

Let $g \in G$ be such that $\rho' = c_{g} \circ \rho$ where $c_{g}: x \mapsto g x g^{-1}$
is the conjugation. Let $\tilde{g} \in \tilde{G}$ be a lift of $g$.
Then both $\tilde{\rho}'$ and $c_{\tilde{g}} \circ \tilde{\rho}$
are lifts of $\rho'$. They must be equal (Proposition ~\ref{T:209}).

\medskip

Step 2: $\tilde{\rho}, \tilde{\rho}'$ are globally conjugate $\quad\Longrightarrow\quad$ $\rho,
\rho'$ are globally conjugate.

\medskip

Easy.

\medskip

Step 3: $\rho, \rho'$ are globally conjugate in image $\quad\Longrightarrow\quad$ $\tilde{\rho},
\tilde{\rho}'$ are globally conjugate in image.

\medskip

Let $g \in G$ be such that $\rho' (H) = g \rho (H) g^{-1}$ and $\tilde{g} \in \tilde{G}$
a lift of $g$.
Then $\tilde{\rho}' (H)$ and $g \tilde{\rho} (H) g^{-1}$ are both connected
subgroups of $\pi_{G}^{-1} (\rho' (H))$. As $\pi_{G}$ is a finite map,
$\tilde{\rho}' (H)$, $g \tilde{\rho} (H) g^{-1}$ and ${(\pi_{G}^{-1} (\rho' (H)))}^{\circ}$
(the identity component of $(\pi_{G}^{-1} (\rho' (H)))$),
must possess the same dimension as $\rho (H)$ and $\rho' (H)$.
Thus these three subgroups must be equal.

\medskip

Step 4: $\tilde{\rho}, \tilde{\rho}'$ are globally conjugate in image $\quad\Longrightarrow\quad$ $\rho,
\rho'$ are globally conjugate in image.

\medskip

Easy also.

\medskip

Step 5: $\rho, \rho'$ are locally conjugate $\quad\Longleftrightarrow\quad$ $\tilde{\rho},
\tilde{\rho}'$ are locally conjugate.

\medskip

Let $T$ be a maximal torus of $H$.

\begin{align}
\text{$\rho, \rho'$ are locally conjugate} &\Longleftrightarrow
\text{$\rho |_{T}, \rho' |_{T}$ are globally conjugate}
&& \text{(Lemma ~\ref{T:211})} \notag \\
&\Longleftrightarrow \text{$\tilde{\rho} |_{T}, \tilde{\rho}' |_{T}$ are globally conjugate}
&& \text{(Step 1 \& 2)} \notag \\
&\Longleftrightarrow \text{$\tilde{\rho}, \tilde{\rho}'$ are locally conjugate}
&& \text{(Lemma ~\ref{T:211})} \notag
\end{align}

\qedsymbol

\medskip

Recall that $M (\rho; G)$ is the set of equivalent classes of $\rho'$:
$H \to G$ which is locally conjugate to $i$, modulo the global conjugacy, and
$M' (\rho; G)$ is the set of equivalent classes of $\rho'$:
$H \to G$ which is locally conjugate to $i$, modulo the global conjugacy in image.
Moreover, $m (\rho; G)$ and $m' (\rho; G)$ as the cardinalities of $M (\rho; G)$
and $M' (\rho; G)$ respectively.

\medskip

\begin{theorem}  \label{T:214}
Let $H$, $G$, $\tilde{H}$, $\tilde{G}$ be connected complex reductive groups,
$\pi_{H}: \tilde{H} \to H$, $\pi_{G}: \tilde{G} \to G$ finite central isogenies,
$\tilde{\rho}$ be a finite central isogenous cover of $\rho$.

\[
\begin{tikzcd}
\tilde{H} \arrow{r}{\rho, (\tilde{\rho}')} \arrow{d}{\pi_{H}}
&\tilde{G} \arrow{d}{\pi_{G}} \\
H \arrow{r}{\rho, (\rho')} &G
\end{tikzcd}
\]

Then the lift of $\rho': H \to G$ to $\tilde{\rho}': \tilde{H} \to \tilde{G}$
define two well defined bijections from $M (\rho, G)$ to $M (\tilde{\rho}, \tilde{G})$
and from $M' (\rho; G)$ to $M' (\tilde{\rho}; \tilde{G})$ respectively.
In particular, $m (\rho; G) = m (\tilde{\rho}; \tilde{G})$ and
$m' (\rho; G) = m' (\tilde{\rho}; \tilde{G})$.
\end{theorem}

\medskip

\emph{Proof of Theorem ~\ref{T:214} and Theorem ~\ref{T:208}}

\medskip

Denote $P (\rho; G)$ the set of homomorphism $\rho': H \to G$
that is locally conjugate to $\rho$ and $P (\tilde{\rho}; \tilde{G})$ similarly.
As $H$ is connected, by Corollary ~\ref{T:212},
we may lifts $\rho'$ to a unique $\tilde{\rho}'$ (see the diagram)
which gives rise to a well defined map
\[
\Phi_{\rho}: P (\rho; G) \to P (\tilde{\rho}; \tilde{G})
\]

\medskip

\emph{Claim}: $\Phi_{\rho}$ is surjective. (In fact, it is bijective.)

\medskip

Granting this claim, and by Theorem ~\ref{T:213}, such $I_{\rho}$ induce
two well defined bijections from $M (\rho, G)$ to $M (\tilde{\rho}, \tilde{G})$
and from $M' (\rho; G)$ to $M' (\tilde{\rho}; \tilde{G})$ with respectively.
Hence we get the theorem.

\medskip

Now we prove the claim. Assume that $\tilde{\rho}'$ is given and it is locally
conjugate to $\tilde{\rho}$. As $\tilde{\rho}$ and $\tilde{\rho}'$ are locally conjugate,
so are $\pi_{G} \circ \tilde{\rho}$ and $\pi_{G} \circ \tilde{\rho}'$, and hence
they should share the same kernel. In particular, both should factor through the
same quotient. Hence the claim and the theorem.

\medskip

Now Theorem ~\ref{T:208}
follows directly with the discussion at the beginning of this subsection.

\medskip

\qedsymbol

\medskip

\begin{lemma}  \label{T:215}
Let $\pi_{G}: \tilde{G} \to G$ be a finite central isogeny
of connected complex reductive groups.
Given $\rho: H \to G$ or $\tilde{\rho}: \tilde{H} \to \tilde{G}$
which is a homomorphism, where $H$ or $\tilde{H}$ is complex connected reductive.
we can complete the following diagram to make $\tilde{\rho}$ a finite central
isogenous cover of $\rho$.
\[
\begin{tikzcd}
\tilde{H} \arrow{r}{\tilde{\rho}} \arrow{d}{\pi_{H}}
&\tilde{G} \arrow{d}{\pi_{G}} \\
H \arrow{r}{\rho} &G
\end{tikzcd}
\]
Moreover, if $G$ and $G'$ be two complex reductive group and $\rho: H \to G$
a homomorphism from a complex connected reductive group $H$
to $G$, then there is a homomorphism $\rho': H' \to G'$
which is isogenous to $\rho$.
\end{lemma}

\medskip

\emph{Proof}: The second assertion follows from the first. Now we work on the first assertion.

\medskip

First, if $\tilde{\rho}$ is given, we can just take $H = \tilde{H}$, $\pi_{H} = {\rm id}$
and $\rho = \pi_{G} \circ \tilde{\rho}$.

\medskip

Next, if $\rho$ is given, we want to get $\tilde{\rho}$.
Let $\tilde{H} = {(H \times_{G} \tilde{G})}^{\circ}$
where $H \times_{G} \tilde{G} = \{(h, \tilde{g}) \in
          G \times \tilde{G} \mid \rho (h) = \pi_{G} (\tilde{g})\}$
and $\tilde{H}$ its identity component,
and $\pi_{H}$ and $\tilde{\rho}$ be the coordinate projection maps
which are homomorphisms.
Then we get the commutative diagram.

\medskip

Now we claim: $\pi_{H}$ is a finite central isogeny. In fact, $\pi_{H}$
is surjective, since $\pi_{G}$ is surjective, and for each $h \in H$,
we can find $\tilde{y}$ with $\pi_{G} (\tilde{y}) = \rho (H)$, and thus
$(h, \tilde{y}) \in \pi_{H}^{-1} (h)$. Moreover, $\pi_{H}$ is a finite central isogeny,
since ${\rm Ker} (\pi_{H}) = 1 \times {\rm Ker} (\pi_{G})$ is finite
and central in $H \times_{G} \tilde{G}$.

\medskip

Done.

\medskip

\qedsymbol

\medskip

\emph{Proof of Theorem ~\ref{TM:A}}:

\medskip

Now let $G' = \mGL (N, \bC), \mSO (2 N + 1, \bC), \mSp (2 N, \bC), \mSO (2 N, \bC)$,
and $G$ an isogenous form of $G'$.
Let $\rho: H \to G$ be a homomorphism of complex connected
reductive groups. Then by Lemma ~\ref{T:215}, there is a $\rho': H' \to G'$,
a homomorphism of complex connected reductive groups, which is isogenous
to $\rho$. By Theorem ~\ref{T:208}, $m (\rho, G) = m (\rho', G')$ and
$m' (\rho, G) = m' (\rho', G')$. Theorem \ref{TM:A} follows now from
Theorem ~\ref{T:202} and Corollary ~\ref{T:204}.

\medskip

\qedsymbol

\bigskip

\section{Main Theorem and Proofs} \label{S:3}

In this section, we mainly study the case $G = \mSO (2 N, \bC)$. In fact, we will
get certain classification of LFMO-special representations.
The definition occurs in both the introduction, and Subsection ~\ref{SS:301}.
In fact, we will prove Theorem ~\ref{TM:B} in Subsection ~\ref{SS:301}
and Theorem ~\ref{TM:C}, Theorem ~\ref{TM:D} and Corollary ~\ref{TM:E}
in Subsection ~\ref{SS:302}.

\medskip

\subsection{$\mSO (2 N)$, Global VS Local, I} \label{SS:301}

% $SO (2 N)$ case, classification, global VS local

\ssmedskip

Now again let $H$ be a complex reductive group,
$G = \mSO (2 N, \bC)$, and $\rho: H \to G$ a homomorphism.

\medskip

Recall some definitions, an element $g \in \mSO (2 N, \bC)$ is said
to \emph{have eigenvalue $1$ or $-1$} if view $g$ as an element
of $\mGL (2 N, \bC)$ via the $2 N$-dimensional standard representation.
$\rho$ (or the representation induced by $\rho$) is said to be \emph{have weight $1$}
if when restricted to one of its maximal torus,
the representation space contains a trivial subrepresentation
as a constituent.

\medskip

\begin{lemma} \label{T:301}
Let $g \in \mSO (2 N, \bC)$ be a semisimple element and $g' = \tau (g)$. Then $g$
and $g'$ are conjugate if and only if $g$ has eigenvalue $\pm 1$.
\end{lemma}

\medskip

\emph{Proof}:

\medskip

Let $V$ be the $2 N$-dimensional complex space with an orthogonal structure $B$
and identify $G = \mSO (2 N, \bC)$ with $\mSO (V, B)$.

\medskip

First we analyze the structure of $V$ under the action of $g$.
As $g$ is semisimple, $V$ is a direct sum of eigenspaces $V_{\lambda}$
of $g$. Moreover, $V_{\lambda}$ and $V_{\lambda'}$ are orthogonal unless
$\lambda = \lambda'^{-1}$, and moreover, $V_{\lambda^{-1}} \cong V_{\lambda}^{*}$
as vector spaces through $B$.

\medskip

Next the only if part. Assume that $g$ and $g'$ are conjugate.
Then $g' = \tau (g)= y g y^{-1}$
for some $y \in \mSO (V, B)$. Thus the automorphism $c_{y^{-1}} \circ \tau$
of $\mSO (V, B)$ fixes $g$. In particular, their is a $x \in \mO (V, B) - \mSO (V, B)$
commutes with $g$. Thus $x$ stabilizes all $V_{\lambda}$.
As $x \in \mO (V, B)$, the action of $x$ on $V_{\lambda}$ and $V_{\lambda^{-1}}$
are contragredient to each other, and hence ${\rm det} (x |_{V_{\lambda^{-1}}}) =
{({\rm det} (x|_{V_{\lambda}}))}^{-1}$.
So $-1 = {\rm det} (x) = {\rm det} (x |_{V_{1}}) {\rm det} (x |_{V_{-1}})$.
So $V_{1}$ or $V_{-1}$ must be nontrivial.

\medskip

Conversely, assume $g$ has an eigenvalue $\epsilon = \pm 1$. Under some basis,
$$g = {\rm Diag} (\alpha_{1}, \ldots, \epsilon,
        \epsilon, \ldots, \alpha_{1}^{-1})$$
and
\[
B = \begin{pmatrix}
0 & 0 & \cdots & 0 & 1 \\
0 & 0 & \cdots & 1 & 0 \\
\vdots & \vdots & \ddots & \vdots & \vdots \\
0 & 1 & \cdots & 0 & 0 \\
1 & 0 & \cdots & 0 & 0
\end{pmatrix}
\]
Then $g$ commutes with the following $z \in \mO (V, B) - \mSO (V, B)$.
\[
z = \begin{pmatrix}
1 & 0 & \cdots & 0 & 0 & 0 & 0 & \cdots & 0 & 0 \\
0 & 1 & \cdots & 0 & 0 & 0 & 0 & \cdots & 0 & 0 \\
\vdots & \vdots & \ddots & \vdots & \vdots & \vdots & \vdots & \ddots & \vdots & \vdots \\
0 & 0 & \cdots & 1 & 0 & 0 & 0 & \cdots & 0 & 0 \\
0 & 0 & \cdots & 0 & 0 & 1 & 0 & \cdots & 0 & 0 \\
0 & 0 & \cdots & 0 & 1 & 0 & 0 & \cdots & 0 & 0 \\
0 & 0 & \cdots & 0 & 0 & 0 & 1 & \cdots & 0 & 0 \\
\vdots & \vdots & \ddots & \vdots & \vdots & \vdots & \vdots & \ddots & \vdots & \vdots \\
0 & 0 & \cdots & 0 & 0 & 0 & 0 & \cdots & 1 & 0 \\
0 & 0 & \cdots & 0 & 0 & 0 & 0 & \cdots & 0 & 1
\end{pmatrix}
\]

\medskip

Now $g' = \tau (g)$ is conjugate to $z g z^{-1} = g$.

\qedsymbol

\medskip

\begin{theorem}  \label{T:302}
Let $H$ be a complex reductive group and
$\rho: H \to G = \mSO (2 N, \bC) \subset \Omega = \mO (2 N, \bC)$
be a homomorphism.
$\rho' = \tau \circ \rho$ for some outer automorphism $\tau$
of $G$. Let $T$ be a maximal torus of $G$. Then

\medskip

\textnormal{(A)} The following are equivalent if $H$ is connected.

\medskip

\textnormal{(a-1)} $\rho$ and $\rho'$ are locally conjugate.

\medskip

\textnormal{(a-2)} $\rho |_{T}$ and $\rho' |_{T}$ are locally conjugate.

\medskip

\textnormal{(a-3)} $C_{\Omega} (\rho (T)) \nsubseteq G$.

\medskip

\textnormal{(a-4)} $\rho$ has weight $1$.

\medskip

\textnormal{(B)} The following are equivalent.

\medskip

\textnormal{(b-1)} $\rho$ and $\rho'$ are globally conjugate.

\medskip

\textnormal{(b-2)} $C_{\Omega} (\rho (H)) \nsubseteq G$.

\medskip

\textnormal{(b-3)} $\rho$ has an irreducible subrepresentation which is both odd dimensional
and of orthogonal type.

\medskip

\textnormal{(C)} The following are equivalent.

\medskip

\textnormal{(c-1)} $\rho$ and $\rho'$ are globally conjugate in image.

\medskip

\textnormal{(c-2)} $N_{\Omega} (\rho (H)) \nsubseteq G$.

\end{theorem}

\medskip

Recall that $C_{G} (G_{0})$ and $N_{G} (G_{0})$ denote the centralizer and the
normalizer of $G_{0}$ in $G$.

\medskip

\emph{Remark}: The equivalence of (a-1), (a-2) and (a-3) apply also
for arbitrary connected group $G$. The equivalence of (b-1) and (b-2),
(c-1) and (c-2) apply also for any group $H$ and $G = \mSO (2 N, \bC)$.

\medskip

\emph{Proof}: Recall that we only assume that $H$ is connected in (A).
Moreover, let $V$ be the representation space of $\rho$ and $B$ the non-degenerate
$H$-invariant bilinear form on $V$ induced by $\rho$.
Let $\Omega = \mO (V, B)$.

\medskip

(A) The equivalence of (a-1), (a-2) follow from Lemma ~\ref{T:211}.
Recall that $\rho' |_{T} = \tau \circ \rho |_{T} = c_{x} \circ \rho |_{T}$
for some $x \in \Omega - G$. Thus (a-3) is equivalent to the fact
that $\rho |_{T} = c_{z} \circ \rho |_{T}$ for some
$z \in \Omega - G$ which is globally conjugate to $\rho' |_{T} = c_{x} \circ \rho |_{T}$
as $x z^{-1} \in G$. Hence this gives rise to (a-1) $\leftrightarrow$ (a-3).

\medskip

Let $t$ be a generic point of $T$.
Assume (a-3), we have $\rho (t)$
commutes with $z \in \Omega - G$,
Hence $\rho (t)$ has eigenvalue $\pm 1$ (Lemma ~\ref{T:301}).
So $V^{T} = V^{\left< t^{2} \right>} \ne 0$ where $\left< t^{2} \right>$
is the group generated by
$t^{2}$ since $t^{2}$ is a also generic point of $T$.
This gives (a-4). So (a-3) implies (a-4).

\medskip

Now assume (a-4), and let again $t$ be a generic point of $T$.
(a-4) implies that $\rho (t)$ has eigenvalue
$1$. Thus $\rho (t)$ commutes with some $z \in \Omega - G$.
As $t$ is a generic point of $T$, $z$ also commutes with
$\rho (T)$, and this gives (a-3). So (a-4) implies (a-3).

\medskip

(B) Recall that $\rho' = \tau \circ \rho = c_{x} \circ \rho$
for some $x \in \Omega - G$. Thus (a-3) is equivalent to the fact
that $\rho = c_{z} \circ \rho$ for some
$z \in \Omega - G$ which is globally conjugate to $\rho' = c_{x} \circ \rho$
as $x z^{-1} \in G$. Hence this gives rise to (b-1) $\leftrightarrow$ (b-2).

\medskip

Now apply our structure theory (Theorem ~\ref{T:205}). For each irreducible
constituent $\sigma$, we have the isotypical component $V_{\sigma}$ of $V$
and $W_{\sigma} = V_{\sigma}$ if $\sigma$ is self-dual and
$V_{\sigma} \oplus V_{\sigma^{\vee}}$ otherwise. Then $V$ is an orthogonal
sum of $W_{\sigma}$. If $\sigma$ is not self-dual, then
$W_{\sigma} = V_{\sigma} \oplus V_{\sigma}$ gives rise to
a complete polarization. If $\sigma$ is self-dual,
$V_{\sigma} \cong W_{\sigma, 0} \otimes W'_{\sigma}$
where $W_{\sigma, 0}$ is equivalent to $\sigma$,
and $H$ acts on the first factor.
Moreover, we have non-degenerate bilinear
forms $B_{\sigma, 0}$, $B'_{\sigma}$ on $W_{\sigma, 0}$
and $W'_{\sigma}$ respectively such that the
isomorphism is an isometry. Moreover
$B_{\sigma, 0}$ and $B'_{\sigma}$ are both symmetric or both
alternating.

\medskip

Assume (b-3), and let $\sigma$ be a irreducible constituent of $V$
which is odd dimensional of orthogonal type. Then $W_{\sigma} = V_{\sigma} \ne 0$
$W'_{\sigma} \ne 0$ is an orthogonal space. Choose an anisotropic vector
$w' \in W'_{\sigma}$ (namely, $B'_{\sigma} (w', w') \ne 0$), and form
$W''$ in $V$ from $W_{\sigma, 0} \otimes w'$ through the isometry. Then $W''$
is a non-degenerate $H$-subspace of $V$ equivalent to $\sigma$. Let $A \in GL (V)$
be such that $A |_{W''} = - {\rm id}$ and $A |_{W''^{\bot}} = 1$.
Then $A$ is $H$-equivariant and $A \in \Omega = \mO (V, B)$
so that $A \in C_{\Omega} (\rho (H))$.
Moreover ${\rm det} (A) = {(-1)}^{{\rm dim} (W'')}
     = {(-1)}^{{\rm dim} \sigma} = -1$ as $\sigma$ is odd dimensional.
Then $A \not\in G$, and this gives (b-2). Hence (b-3) implies (b-2).

\medskip

Now assume that (b-3) fails, namely all irreducible constituents of $\rho$
are either not self-dual or even dimensional. We want to
prove $C_{\Omega} (\rho (H)) \subset G$ so that (b-2) fails.

\medskip

Choose arbitrary $T \in C_{\Omega} (\rho (H))$. Then $T$ is $H$-equivariant, and hence
$T$ stabilizes all $V_{\sigma}$ and $W_{\sigma}$. Denote $T_{\sigma} = T |_{W_{\sigma}}$.
We \emph{claim} that ${\rm det} (T_{\sigma}) = 1$. Once we have the claim,
$T \in \Omega \cap \mSL (V) = G$. Hence $C_{\Omega} (\rho (H)) \subset G$.

\medskip

Assume first that $\sigma$ is not self-dual.
Then $T_{\sigma}$ stabilizes $V_{\sigma}$ and $V_{\sigma^{\vee}}$, which give rise to
a complete polarization of $W_{\sigma}$.
Through $B$, $V_{\sigma^{\vee}} \cong V_{\sigma}^{*}$. As $T$ is
orthogonal, then under a choice of dual bases of $V_{\sigma}$ and
$V_{\sigma^{\vee}}$, the matrix representations of $T$
on these two spaces are transpose inverse to each other. Hence
${\rm det} (T_{\sigma}) = 1$.

\medskip

Assume now that $\sigma$ is self-dual. Then $\sigma$ must be even dimensional
and $V_{\sigma} = W_{\sigma} \cong W_{\sigma, 0} \otimes W'_{\sigma}$.
As $T$ is $H$-equivariant, then by Schur's lemma (Lemma ~\ref{T:206} also,
plus a good exercise in linear algebra),
the action of $T$ on $V_{\sigma}$ is induced by $1 \otimes T'_{\sigma}$
where $T'_{\sigma} \in \mO (W'_{\sigma}, B'_{\sigma})$.
In particular ${\rm dim} (T'_{\sigma}) = \pm 1$ and
\[{\rm dim} (T_{\sigma}) =
{({\rm dim} (T'_{\sigma}))}^{{\rm dim} (W_{\sigma, 0})} = {(\pm 1)}^{{\rm dim} (\sigma)} = 1\]

\medskip

In all cases, ${\rm det} (T_{\sigma}) = 1$. So ${\rm det} (T) = 1$
and $T \in \Omega \cap \mSL (V) = G$. Thus (b-2) fails. So
(b-2) implies (b-3).

\medskip

(C) $\rho$ and $\rho' = \tau \circ \rho$ are globally conjugate in image in $G$
if and only if $\rho (H)$ and $\tau (\rho (H))$ are conjugate,
if and only if $\rho (H)$ is stabilized by $c_{y} \circ \tau$ for some $y \in G$.
As $\tau = c_{x}$ for some $x \in \Omega - G$, the last situation
occurs if and only if $H$ is normalized by some element
in $\Omega - G$. So (c-1) and (c-2) are equivalent.

\medskip

\qedsymbol

\medskip

\emph{Proof of Theorem ~\ref{TM:B}}:

\medskip

Now (a) is Theorem ~\ref{T:302} (A) (a-1) $\leftrightarrow$ (a-4),
(b) is Theorem ~\ref{T:302} (B) (b-1) $\leftrightarrow$ (b-3),
(c) is Theorem ~\ref{T:302} (A) (c-1) $\leftrightarrow$ (c-2).

\medskip

\qedsymbol

\medskip

\begin{corollary}  \label{T:303}
Let $H$ be ac complex reductive group and
$\rho: H \to G = \mSO (2 N, \bC)$ be a homomorphism.
Then $m (\rho; G), m' (\rho; G) = 1, 2$. Moreover, $m (\rho; G) = 2$
if and only if \textnormal{(A)} holds and \textnormal{(B)} fails,
and $m' (\rho; G) = 2$ if and only if \textnormal{(A)} holds
and \textnormal{(C)} fails.
Here \textnormal{(A)}, \textnormal{(B)} and \textnormal{(C)}
are the any one of equivalent conditions of
\textnormal{(A)}, \textnormal{(B)} and \textnormal{(C)} in Theorem ~\ref{T:302}.
In particular, if $\rho$ is stable, then $m (\rho; G) = 2$ if and only if
$\rho$ has weight $1$.
\end{corollary}

\medskip

\emph{Proof}: The first condition is exactly Corollary ~\ref{T:204},
and the rest can be routinely checked easily.

\qedsymbol

\bigskip

\subsection{$\mSO (2 N)$, Local VS Global, II} \label{SS:302}

\ssmedskip

In this part, we focus on LFMO-special representation. Let $H$
be a complex (reductive) group.
Recall that a homomorphism $\rho: H \to \mO (N, \bC), \mSp (2 N, \bC)$
(or the induced self-dual representation) is said to be stable if the representation it
induced is irreducible, and is said to be essential if it has no
nonzero totally isotropic
constituent. For definition, see the introduction, or Subsection ~\ref{SS:202}
(after Corollary ~\ref{T:203} and before the proof of Theorem ~\ref{T:202}).

\medskip

% $SO (2 N)$ case, classification, global VS local, in image I.

Recall the definition. A homomorphism $\rho: H \to \mSO (2 N, \bC)$
is said to be \emph{LFMO-special} if, for
 $\rho' = \tau \circ \rho$ for an outer automorphism $\tau$ of $\mSO (2 N)$,
the following conditions for $\rho$ and $\rho'$ hold.

\medskip

(1) $\rho$ is \emph{essential}, namely, the image of $\rho$ is not contained in any parabolic
subgroup of $\mSO (2 N)$.

\medskip

(2) $\rho$ and $\rho'$ are locally conjugate.

\medskip

(3) $\rho$ and $\rho'$ are not globally conjugate in image.

\medskip

In this case, we also say that the induced $2 N$-dimensional representation
is LFMO-special.

\medskip

The name ``LFMO'' comes from the following:
If $\rho$ is LFMO-special, then $\rho$ and $\rho'$ will lead to failure of
multiplicity one. That is also why we call this name (see \cite{Wang2011}).
In one word, $\rho$ is LFMO-special if and only if
$\rho$ is essential and $m' (\rho; \mSO (2 N, \bC)) = 2$.

\bigskip

\emph{Proof of Theorem ~\ref{TM:C}}: Easy right now.

\medskip

According to our definition, the condition (1), (2) and (3) of Theorem ~\ref{TM:C}
match the condition (1), (2) and (3) of our definition
(where (1) from Proposition ~\ref{T:207},
(2) and (3) from Theorem ~\ref{T:302}, also see Corollary ~\ref{T:303}).
Moreover, the last statement follows from Theorem ~\ref{T:208}
as when view $\rho, \rho'$ as homomorphisms from $H$ to $G$,
they are isogenous since they are quasi-equivalent when
viewed as complex representations.

\qedsymbol

\medskip

% $SO (2 N)$ case, classification, stable case, in image II.

Now we focus on stable orthogonal representations and try to
classify stable LFMO-special representations.
Note that from Corollary ~\ref{T:303}, we have $m (\rho; G) = 2$
if and only if $\rho$ has weight $1$,
and now $m' (\rho; G) = 2$ if and only if $\rho$ is LFMO-special.
As each $H$ has a finite central isogenous cover
$T \times H_{1} \times H_{2} \times \ldots \times H_{r}$ which is a decomposable
group for $T$ a torus, and $H_{i}$ simple Lie group for each $i$, $\rho$
is isogenous to $\rho'$ a homomorphism from $H'$ to $G$.
Theorem ~\ref{TM:C}, Theorem ~\ref{T:208} and Corollary ~\ref{T:303}
enable us to reduce the classification problem to the case when
$H$ itself is decomposable. Moreover, $\rho$ is also factorizable
as a tensor product of self-dual representations, and the factorization is unique.
and $\rho |_{T}$ is trivial as $T$ is central and $\rho$ is stable (Lemma ~\ref{T:206}).
We will observe this scenario (see Assumption A below),
and replace each $H_{i}$ by its isogenous form when convenient
till the end of this section.

\medskip

\emph{Remark}: When $H$ is decomposable, $H$ is semisimple if and only if $T = 1$;
In addition, $H$ is simply connected / adjoint if and only if so are $H_{i}$
for all $i$.

\medskip

\begin{lemma} \label{T:304}
With $H$ decomposable as above, $\rho$ has weight $1$ if and only if so does $\rho_{i}$
for each $i \ge 1$.
\end{lemma}

\medskip

\qedsymbol

\medskip

The lemma above plus Theorem ~\ref{TM:C}, enable us to reduce the problem
to the study of the lift from automorphisms from $\rho (H)$
to the automorphism of $G$.

\medskip

\begin{lemma}  \label{T:305}
Let $H$ be a complex semisimple, simply connected or adjoint group,
and $\rho: H \to G = \mSp (2 N, \bC) \,\text{or}\, \mSO (N, \bC)$
then automorphisms of $\rho (H)$ lift to those of $H$, i.e.,
given $\phi \in {\rm Aut} (\rho (H))$, there exists a $\phi_{0} \in {\rm Aut} (H)$,
such that $\phi \circ \rho = \rho \circ \phi_{0}$.
\[
\begin{tikzcd}
H \arrow{r}{\rho} \arrow{d}{\phi_{0}}
&\rho (H) \arrow{d}{\phi} \\
H \arrow{r}{\rho} &\rho (H)
\end{tikzcd}
\]
\end{lemma}

\medskip

This lemma enable us to reduce our problem to the study of ${\rm Aut} (H)$
in place of ${\rm Aut} (\rho (H))$.

\medskip

\emph{Proof}: First case, assume that ${\rm Ker} (\rho)$ is finite. $\rho$
is also a finite central isogeny from $H$ to $\rho (H)$. Thus
we get $\phi_{0}: H \to H$ as a lift
of $\phi \circ \rho: H \to \rho (H)$, since when $H$ is simply connected
we can do this definitely (see the proof of Proposition ~\ref{T:209}),
and when $H$ is adjoint, $\rho$ is an isomorphism.

\medskip

In general, let $N_{0} = {({\rm Ker} (\rho))}^{\circ}$.
Then $N_{0}$ is normal in $H$. As $H$ is semisimple and simply connected / adjoint,
$H = N_{0} \times H'$ for some $H'$ which is also normal and semisimple.
In fact both $N_{0}$ and $H'$ are product of simply connected / adjoint simple groups.
Now $\rho$ factors through $H' \cong H / N_{0}$ and
induces a finite central isogeny from $H'$ to $\rho (H)$.
Thus $\phi \in {\rm Aut} (\rho (H))$ lifts to $\phi'_{0} \in {\rm Aut} (H')$
from the first case. The lemmas follows now easily when take
$\phi_{0} = {\rm id_{N_{0}}} \times \phi'_{0}$.

\medskip

\qedsymbol

\medskip

Before we state our main results, we list some definition (which is also used
in \cite{Wang2007}, \cite{Wang2011}).

\medskip

Let $\rho_{0}: H \to \mGL (V)$ be a finite dimensional representation,
and $V$ its representation space.
We say that $\phi \in {\rm Aut} (H)$ is \emph{$\rho_{0}$-liftable}
if $\rho_{0}$ and $\rho_{0} \circ \phi$ are equivalent,
i.e., there is a $T \in \mGL (V)$
such that $\rho_{0} \circ \phi = c_{T} \circ \rho_{0}$, namely,
\[
\begin{tikzcd}
H \arrow{r}{\rho_{0}} \arrow{d}{\phi}
&\mGL (V) \arrow{d}{c_{T}} \\
H \arrow{r}{\rho_{0}} &\mGL (V)
\end{tikzcd}
\]
and in this case we also say that \emph{$\phi$ lifts to $T$ (through $\rho$)}.

\medskip

Let $\rho: H \to G = \mSO (N, \bC)$ or $\mSp (2 N, \bC) < \mGL (V)$
be a homomorphism with $V$ its representation space, and also denote
$\rho$ as the representation.
We say that $\phi$ is \emph{$\rho$-even}
if $G = \mSO (2 N, \bC)$ or $\mSp (2 N, \bC)$ and $\phi$ lifts
to some $T \in G$,
$\phi$ is \emph{$\rho$-odd} if $G = \mSO (2 N, \bC)$, and
$\phi$ lifts to some $T \in \Omega - G$ where $\Omega = \mO (2 N, \bC) < \mGL (V)$,
and $\phi$ is \emph{$\rho$-neutral} if $G = \mSO (2 N + 1, \bC)$ and
$\phi$ is $\rho$-liftable.

\medskip

Also, recall that two representations
$\rho$ and $\rho'$ of $H$ and $H'$
are \emph{quasi-equivalent} if and only if there are
finite central isogenies $\iota: H'' \to H$ and
$\iota: H'' \to H'$ such that as representations of $H''$,
$\rho \circ \iota$ and $\rho' \circ \iota'$
are conjugate.

\medskip

\begin{proposition} \label{T:306}
Let $H$ be a complex connected reductive group,
$\rho: H \to G = \mSp (2 N, \bC)$ or $\mSO (N', \bC)$
be a stable homomorphism and
$\phi \in {\rm Aut} (H)$. Assume that $\phi$ is $\rho$-liftable.
Then $\phi$ is either $\rho$-even, or $\rho$-odd, or
$\rho$-neutral, and only one case occurs.
Moreover, if $H$ is semisimple and simply connected (or adjoint) and $G = \mSO (2 N, \bC)$,
then $\rho$ is LFMO-special
if and only if: \textnormal(1) $\rho$ has weight $1$ and \textnormal{(2)} each
$\phi \in {\rm Aut} (H)$ is either $\rho$-even or not $\rho$-liftable.
\end{proposition}

\medskip

\emph{Proof}:

\medskip

Let $B$ be a non-degenerate $H$-invariant symmetric / alternating bilinear form on
$V$, the representation space of $\rho$, and $\Omega = \mO (V, B)$.
Assume that $\phi \in {\rm Aut} (H)$ is $\rho$-liftable, and $\phi$
lifts to $T \in \mGL (V)$. As $\rho$ is stable and $V$ is irreducible,
$B^{T} = C B$ for some $C \in \bC$
where $B^{T} (u, v) = B (T u, T v)$
(Schur's lemma, Lemma ~\ref{T:206}). Thus $\phi$ lifts to
$C^{-1} T \in \mO (V, B)$. If $B$ is alternating, then $\Omega = G = \mSp (2 N, \bC)$
and $\phi$ is $\rho$-even; If $B$ is symmetric and $V$ is odd dimensional,
then $\Omega = G = \mSO (2 N + 1, \bC)$, and $\phi$ is $\rho$-neutral;
If $B$ is symmetric and $V$ is even dimensional,
then $\Omega = \mO (2 N, \bC)$, and $\phi$ is either $\rho$-even
or $\rho$-odd. Moreover, if $\phi$ lifts to $T, T' \in \Omega$,
then $T^{-1} T'$ centralizes $\rho (H)$. As $V$ is irreducible,
$T^{-1} T' \in \Omega$ is a scalar,
and hence is $\pm 1 \in G = \mSO (2 N, \bC)$.
Thus $T \in G$ if and only if $T' \in G$. So $\phi$
can't be $\rho$-even or $\rho$-odd at the same time.

\medskip

Next assertion: Assume now that $\rho$ has weight $1$. Since $\rho$
is assumed to be stable, and hence essential, by Theorem ~\ref{TM:C},
it suffices for us to show that the condition (3): If an automorphism
$\phi'$ of $\rho (H)$ lifts to $T \in \Omega$ then $T \subset G$
is equivalent to the following:

\medskip

Condition (x): all $\rho$-liftable automorphisms of
$H$ are $\rho$-even.

\medskip

First (3) implies (x). Assume that $\phi \in {\rm Aut} (H)$ is $\rho$-liftable,
and lifts to $T \in \Omega$. Then the conjugation $c_{T}$ induces
an automorphism of $\rho (H)$, and thus by (3), $c_{T}$ must be inner,
and hence $\phi$ is $\rho$-even.

\medskip

Conversely, (x) implies (3). Assume that $\phi' \in {\rm Aut} (\rho (H))$
lifts to an automorphism $L$ on $G$. Then by Lemma ~\ref{T:305}
and the assumption that $H$ is simply connected or adjoint,
there is a $\phi \in {\rm Aut} (H)$ such that
$\rho \circ \phi = \phi' \circ \rho$ (also see the diagram in the proof
of Lemma ~\ref{T:305}). Thus $\phi$ is $\rho$-liftable, and by (x),
is $\rho$-even. and hence $L$ must be inner, and hence a conjugation
on $G$.

\medskip

\qedsymbol

\medskip

\begin{lemma}  \label{T:307}
Let $H$ be a semisimple, simply connected or adjoint Lie group,
and $H = H_{1} \times H_{2} \times \ldots \times H_{r}$ with $H_{i}$
simple. Then ${\rm Aut} (H)$ is generate by the following families:

\medskip

\textnormal{(Type 1: Decomposable)}: $\phi_{1}
     \times \phi_{2} \times \ldots \times \phi_{r}$
where $\phi_{i} \in {\rm Aut} (H_{i})$.

\medskip

\textnormal{(Type 2: Swapping)}:
\begin{align}
T_{\lambda} &(g_{1}, g_{2}, \ldots, g_{i}, \ldots, g_{j}, \ldots, g_{r})
\notag \\
&= (g_{1}, g_{2}, \ldots, \lambda^{-1} (g_{j}), \ldots, \lambda (g_{i}), \ldots, g_{r})
\notag
\end{align}
where $\lambda: H_{i} \overset{\cong}{\to} H_{j}$ is an isomorphism.

\end{lemma}

\medskip

\emph{Proof}: Let $T \in {\rm Aut} (G)$. Then $T (H_{i})$ is a normal simple
subgroup of $H$. As $H$ is simply connected /adjoint, $T (H_{i}) = H_{j}$
for some $j$.
Put $\sigma \in S_{n}$ such that $T (H_{i}) = H_{\sigma (i)}$. Then
$H_{i} \cong H_{\sigma (i)}$. Thus for each cycle $C$ in $\sigma$
and $H_{i}$ are isomorphic all indices $i$ occurred $C$.
Thus, there is a product of Type 2 automorphisms
$T'$ with $T' (H_{i}) = T (H_{i})$ for all $i$.
Note that $T'^{-1} \circ T$ is of Type 1.

\medskip

\qedsymbol

\medskip

\begin{theorem}  \label{T:308}
Let $H$ be a semisimple, simply connected or adjoint Lie group,
and $H = H_{1} \times H_{2} \times \ldots \times H_{r}$ with $H_{i}$
simple. Let $\rho = \bigotimes_{i} \rho_{i}$ where $\rho_{i}$
is a finite dimensional irreducible complex representation of $H_{i}$.
Then ${\rm Aut} (H, \rho)$, the set of $\rho$-liftable
automorphisms of $H$, is generated by the following families

\medskip

\textnormal{(Type 1L: Decomposable)}: $\phi_{1}
     \times \phi_{2} \times \ldots \times \phi_{r}$
where $\phi_{i} \in {\rm Aut} (H_{i}, \rho_{i})$.

\medskip

\textnormal{(Type 2L: Swapping)}:
\begin{align}
T_{\lambda} &(g_{1}, g_{2}, \ldots, g_{i}, \ldots, g_{j}, \ldots, g_{r})
\notag \\
&= (g_{1}, g_{2}, \ldots, \lambda^{-1} (g_{j}), \ldots, \lambda (g_{i}), \ldots, g_{r})
\notag
\end{align}
where $\lambda: H_{i} \overset{\cong}{\to} H_{j}$ is any isomorphism,
and $\rho_{i}, \rho_{j} \circ \lambda$ are equivalent.

\end{theorem}

\medskip

\begin{lemma}  \label{T:309}
Let $H$ and $H'$ be two semisimple Lie groups,
which are both simply connected or adjoint,
$\rho, \rho'$ finite dimensional complex
representations of $H$ and $H'$ respectively.
Then the following are equivalent: \textnormal{(a)} $\rho$ and $\rho'$
are quasi-equivalent. \textnormal{(b)} There is an isomorphism $\lambda: H \to H'$
such that $\rho$ and $\rho' \circ \lambda$ are equivalent.
\end{lemma}

\medskip

\emph{Proof}: It is obvious that (b) implies (a). Now assume (a). Then
$H$ and $H'$ are isogenous and we have the following diagram
\[
\begin{tikzcd}
H'' \arrow{r}{\pi} \arrow{dr}{\pi'}
& H \arrow{r}{\rho}
&\mGL (V) \arrow{d}{\cong} \\
& H' \arrow{r}{\rho'} &\mGL (V)
\end{tikzcd}
\]
for some finite central isogenies $\pi$ and $\pi'$. If both $H$ and $H'$
are simply connected, then $\pi$ and $\pi'$ are isomorphisms. If both
$H$ and $H'$ are adjoint, then by factoring by the center,
can choose $H''$ also adjoint, so that
$\pi$ and $\pi'$ are isomorphisms.

\medskip

\qedsymbol

\medskip

\begin{proposition} \label{T:310}
Let $H, H'$ be two semisimple complex Lie groups and $\rho,
\rho'$ two finite dimensional irreducible complex representations of $H, H'$
respectively. Let $\phi \in {\rm Aut} (H)$ and $\phi' \in {\rm Aut} (H')$.

\medskip

\textnormal{(1)} $\phi \times \phi'$ is $\rho \otimes \rho'$-liftable
if and only if
$\phi$ is $\rho$-liftable and $\phi'$ is $\rho'$-liftable.

\medskip

\textnormal{(2)} Assume that $\lambda: H \to H'$ is an isomorphism.
Then $T_{\lambda}: H \times H' \to H \times H'$ defined as
\[
T_{\lambda} (h, h') = (\lambda^{-1} (h'), \lambda (h))
\]
is $\rho \otimes \rho'$-liftable if and only if $\rho,
\rho' \circ \lambda$ are equivalent. In particular, if $\rho$
and $\rho'$ are quasi-equivalent, and both $H$
and $H'$ are simply connected or adjoint, then there is an isomorphism $\lambda:
H \to H'$ such that $\rho = \rho' \circ \lambda$ and $T_{\lambda}$
is $\rho \otimes \rho'$-liftable.

\medskip

Now in addition, assume that $\rho$ and $\rho'$ are self-dual.
We have:

\medskip

\textnormal{(3)} Assuming that $\phi'$ is $\rho'$-neutral.
Then $\phi \times \phi'$ is $\rho \otimes \rho'$-neutral
(resp.\ $\rho \otimes \rho'$-even, / $\rho \otimes \rho'$-odd)
if and only if
$\phi$ is $\rho$-neutral (resp.\ $\rho$-even, $\rho$-odd).

\medskip

\textnormal{(4)} If $\phi$ are $\rho$-odd or $\rho$-even and
$\phi'$ are $\rho'$-odd or $\rho'$-even, then $\phi \times \phi'$
are $\rho \otimes \rho'$-even.

\medskip

\textnormal{(5)} Same assumption as in \textnormal{(2)}.
Then $T_{\lambda} (h, h')$ is $\rho \otimes \rho'$-neutral
if and only if $\rho$ is odd dimensional, $\rho \otimes \rho'$-even
if and only if ${\rm dim} (\rho) \equiv 0 \pmod{4}$ and
$\rho \otimes \rho'$-odd
if and only if ${\rm dim} (\rho) \equiv 2 \pmod{4}$.

\end{proposition}

\medskip

\emph{Proof}:

\medskip

(1) $\phi \times \phi'$ is $\rho \otimes \rho'$-liftable
if and only if
$(\rho \otimes \rho') \circ (\phi \times \phi')
= (\rho \circ \phi) \otimes (\rho' \circ \phi')$
and $\rho \otimes \rho'$ are equivalent,
if and only if
$\rho \circ \phi$ and $\rho$ are equivalent and
$\rho' \circ \phi'$ and $\rho'$ are equivalent,
if and only if $\phi$ is $\rho$-liftable and $\phi'$ is $\rho'$-liftable.

\medskip

(2) First, assume that $T_{\lambda}$ is $\rho$-liftable.
Then $(\rho \otimes \rho') \times T_{\lambda}
    = (\rho' \circ \lambda^{-1}) \otimes (\rho \circ \lambda)$
and $\rho \otimes \rho'$ are equivalent, and
then $\rho$ and $\rho' \circ \lambda$ are then equivalent.

\medskip

Next, assume that $\rho$ and $\rho' \circ \lambda$ are equivalent.
Then $\rho \circ \lambda^{-1}$ and $\rho'$
are also equivalent, and hence $(\rho \otimes \rho') \times T_{\lambda}
    = (\rho' \circ \lambda^{-1}) \otimes (\rho \circ \lambda)$
and $\rho \otimes \rho'$ are equivalent.
Hence $T_{\lambda}$ is $\rho \times \rho'$-liftable.

\medskip

Finally, if $\rho$ and $\rho'$ are quasi-equivalent and
$H$ and $H'$ are both simply connected or both adjoint,
there is a $\lambda: H \overset{\cong}{\to} H'$
such that $\rho = \rho' \circ \lambda$, and hence $T_{\lambda}$
is $\rho \otimes \rho'$-liftable.

\medskip

(3) \& (4) Let $V, V'$ be the representation spaces of $\rho, \rho'$
respectively, and $B$ the non-degenerate $H$-invariant bilinear form of $V$
and $B'$ the non-degenerate $H'$-invariant bilinear form of $V'$.
Let $\Omega = \mO (V, B)$ and $\Omega' = \mO (V', B')$.
Assume that $\phi, \phi'$ lift to $A \in \Omega$,
$A' \in \Omega'$ (see Proposition ~\ref{T:306}).

\medskip

Then $\phi \times \phi'$ lift to $A \otimes A' \in \tilde{\Omega}
 = \mO (V \otimes V', B \otimes B')$.
Note that
\[
{\rm det} (A \otimes A') = {{\rm det} (A)}^{{\rm dim} (\rho')}
{{\rm det} (A')}^{{\rm dim} (\rho)}
\]
while ${\rm det} (A)$ and ${\rm det} (A')$ are $\pm 1$.

\medskip

Now we divide in cases.

\medskip

(Case 3-1): Both $\rho$ and $\rho'$ are odd dimensional.
In this case, ${\rm dim} (\rho \otimes \rho')$ is odd,
and hence if $\phi$ is $\rho$-neutral, then $\phi \times \phi'$ is
$\rho \otimes \rho'$-neutral.

\medskip

(Case 3-2): $\rho'$ is odd dimensional and $\rho'$ is symplectic.
In this case, $\rho \otimes \rho'$ is also symplectic,
and hence $\phi$ is $\rho$-even and $\phi \times \phi'$ is
$\rho \otimes \rho'$-even.

\medskip

(Case 3-3): $\rho'$ is odd dimensional and $\rho'$
is even dimensional and orthogonal.
In this case, $\rho \otimes \rho'$ is also
even dimensional orthogonal. Moreover,
${\rm det} (A \otimes A') = {\rm dim} (A)$.
as $\phi \times \phi'$ lifts to $A \otimes A' \in \tilde{\Omega}
    = \mO (V \otimes V', B \otimes B')$.
Hence if $\phi$ is $\rho$-even (resp.\ $\rho$-odd),
then $A \in SO (V, B)$,
${\rm det} (A)$ is $1$ (resp.\ $-1$),
${\rm det} (A \otimes A')$ is $1$ (resp.\ $-1$),
$A \otimes A' \in \mSO (V \otimes V', B \otimes B')$,
and thus $\phi \times \phi'$ is $\rho \otimes \rho'$-even.

\medskip

(Case 4): All $\rho$ and $\rho'$ are even dimensional.
In this case, $\rho \otimes \rho'$ is also even dimensional.
Moreover,
$\phi \times \phi'$ lifts to $A \otimes A' \in \tilde{\Omega}
 = \mO (V \otimes V', B \otimes B')$
and ${\rm det} (A \otimes A') = 1$.
Hence $A \otimes A' \in \mSO (V \otimes V', B \otimes B')$,
and thus $\phi \times \phi'$ is $\rho \otimes \rho'$-even.

\medskip

(5) Now by assumptions in (2), $\rho$ and $\rho' \circ \lambda$ are conjugate,
where $\lambda : H \overset{\cong}{\to} H'$.
Hence we have an $H$-equivariant isometry $E$ from $V$ to $V'$
(Proposition ~\ref{T:202}).

\medskip

Then by (2),
$T_{\lambda}$ is $\rho \otimes \rho'$ liftable, and it
lifts to some isometry $A_{\lambda}: V \otimes V' \to V \otimes V'$,
where $A_{\lambda} (v \otimes v') = (E^{-1} (v') \otimes E (v))$
for each $v \in V, v' \in V'$. As ${\rm det} (A_{\lambda}) = {(-1)}^{n (n - 1)}$
where $n = {\rm dim} (\rho)$, then it is $-1$ if $n \equiv 2 \pmod{4}$, and
$-1$ if $n \equiv 0 \pmod{4}$. (Recall that the signature of $\sigma: X \times X \to X \to X$
that sends $(x, y)$ to $(y, x)$ is ${(-1)}^{|X| (|X| - 1)}$.)
Thus, $T_{\lambda}$ is $\rho \otimes \rho'$-neutral
if ${\rm dim} (\rho)$ is odd,
$T_{\lambda}$ is $\rho \otimes \rho'$-odd
if ${\rm dim} (\rho) \equiv 2 \pmod{4}$, and
$T_{\lambda}$ is $\rho \otimes \rho'$-even
if ${\rm dim} (\rho) \equiv 0 \pmod{4}$.

\medskip

\qedsymbol \medskip

\emph{Proof of Theorem ~\ref{T:308}}: We proceed in steps.

\medskip

(Step 1): In two families in Lemma ~\ref{T:307},
the $\rho$-liftable ones are exactly in Type 1L and Type 2L.

\medskip

First, Let $\phi = \phi_{1} \times \ldots \times \phi_{r} \in {\rm Aut} (H)$
of Type 1. Then by Proposition ~\ref{T:310} (1),
$\phi$ is $\rho$-liftable if and only if
$\phi_{i}$ is $\rho_{i}$-liftable.
So Type 1 $\&$ $\rho$-liftable $\Leftrightarrow$ Type 1L.

\medskip

Next, like the proof of Proposition ~\ref{T:310} (2),
let $\phi = T_{\lambda}$ of Type 2 as in Lemma ~\ref{T:307}, and
$\rho' = \rho \circ \phi$. Then $\rho' = \bigotimes_{i} \rho'_{i}$
where $\rho'_{i} = \rho_{j} \circ \lambda$, $\rho'_{j} = \rho_{i} \circ \lambda^{-1}$,
and $\rho'_{k} = \rho_{k}$ for $k \ne i, j$.
$\phi$ is $\rho$-liftable, if and only if $\rho$ and $\rho'$ are equivalent,
if and only if $\rho_{k}$ and $\rho'_{k}$ are equivalent,
if and only if $\rho_{i}$ and $\rho_{j} \circ \lambda$, $\rho_{j}$
and $\rho_{i} \circ \lambda^{-1}$ are equivalent, if and only if
$\rho_{i}$ and $\rho_{j} \circ \lambda$ are equivalent.
Hence Type 2 $\&$ $\rho$-liftable $\Leftrightarrow$ Type 2L.

\medskip

(Step 2): General case. For $\phi \in {\rm Aut} (H)$,
we have $\sigma \in S_{n}$ such that $\phi (H_{i}) = H_{\sigma (i)}$
(see Lemma ~\ref{T:307} and its proof).
Thus by the same argument as last step,
$\rho_{i}$ and $\rho_{j} \circ \phi$ are equivalent, and hence
$H_{i}$ and $H_{j}$ are quasi-equivalent. Thus
for all $i$ occurred in a cycle $C$ of $\sigma$, $H_{i}$
are isomorphic since they are isogenous and both simply connected or both
adjoint. Write $\sigma = t_{1} t_{2} \ldots t_{l}$ for transposition $t_{i}$
where $i_{1}, i_{2}$, all indices of each $t_{i} = (i_{1} i_{2})$
occur in the same cycle in $\sigma$.
Then from Step 1, there is a $\phi_{i} = T_{\lambda_{i}}$ of Type 2L, where
$\lambda_{i}: H_{i_{1}} \to H_{i_{2}}$ an isomorphism.
Let $T = \phi_{l} \circ \phi_{l - 1} \circ \ldots \circ \phi_{1}$.
Then $T$ is $\rho$-liftable and $T (H_{i}) = \phi (H_{i})$ for all $i$.
Thus $T^{-1} \phi$ is of Type 1, and hence Type 1L. Hence the theorem.

\medskip

\qedsymbol

\medskip

\begin{proposition} \label{T:311}
All notations and assumptions as in Theorem ~\ref{T:308}. Assume that
$\rho$ are self-dual. Then

\medskip

\textnormal{(Type 1L)}: Let $\phi = \phi_{1} \times \ldots \times \phi_{r}$
and assume that it is $\rho$-liftable.
Then $\phi$ is $\rho$-neutral if and only if $\rho_{i}$ is odd dimensional
for all $i$. $\phi$ is $\rho$-odd if and only if for exactly one $i$,
$\phi_{i}$ are $\rho_{i}$-odd, and for all other $j$, $\rho_{j}$ are odd dimensional.
$\phi$ is $\rho$-even if and only one of the following happens:
Either for some $i$, $\phi_{i}$ is $\rho_{i}$-even, or for at least two
$i \ne j$, $\rho_{i}$ and $\rho_{j}$ are even dimensional.

\medskip

\textnormal{(Type 2L)}: Let
\begin{align}
T_{\lambda} &(g_{1}, g_{2}, \ldots, g_{i}, \ldots, g_{j}, \ldots, g_{r})
\notag \\
&= (g_{1}, g_{2}, \ldots, \lambda^{-1} (g_{j}), \ldots, \lambda (g_{i}), \ldots, g_{r})
\notag
\end{align}
where $\lambda: H_{i} \overset{\cong}{\to} H_{j}$ is any isomorphism.
Assume that $\rho_{i}, \rho_{j} \circ \lambda$ are equivalent.
Then $T_{\lambda}$ is $\rho$-neutral if and only if for all $k$, $\rho_{k}$
are odd dimensional, $T_{\lambda}$ is $\rho$-odd if and only if
${\rm dim} (\rho_{i}) \equiv 2 \pmod{4}$, and for all $k \ne i, j$,
$\rho_{k}$ is odd dimensional, $T_{\lambda}$ is $\rho$-even if and only if
${\rm dim} (\rho_{i}) \equiv 0 \pmod{4}$, or for some $k \ne i, j$,
$\rho_{k}$ is even dimensional.
\end{proposition}

\medskip

\emph{Proof}:

\medskip

(1) First Type 1L. This follows from Proposition ~\ref{T:310} (1), (3) and (4).

\medskip

(2) Now work for Type 2L. If $r = 2$, then the assertions follow from
Proposition ~\ref{T:310} (2), (5). For general $r$, without loss of
generality may assume that $i = 1, j = 2$.
Let $H' = H_{1} \times H_{2}$, $H'' = H_{3} \times \ldots \times H_{r}$,
$\rho' = \rho_{1} \otimes \rho_{2}$,
$H'' = \rho_{3} \otimes \ldots \otimes \rho_{r}$
and $\phi' = T_{\lambda} |_{H'}$, and $\phi = T_{\lambda} = \phi' \times 1$
on $H = H' \times H''$. Then $\phi'$ is $\rho'$-liftable.
Now Work in case by case carefully. $\phi$ is $\rho$-neutral
if and only if all $\rho_{k}$ are odd dimensional ((1) in this proof).
If $\rho''$ is odd dimensional, then $\phi$ is $\rho$-odd (resp.\ $\rho$-even)
if and only if $\phi'$ is $\rho'$-odd (resp.\ $\rho'$-even) ((1) in this proof),
if and only if ${\rm dim} (\rho_{1}) \equiv 2 \pmod{4}$
(resp.\ ${\rm dim} (\rho_{1}) \equiv 0 \pmod{4}$).
If $\rho''$ is even dimensional, then $1|_{H''}$ is always
$\rho''$-even, $\phi$ is then $\rho$-even.

\medskip

\qedsymbol

\medskip

\begin{theorem}  \label{T:312}
Let $H$ be a semisimple, simply connected or adjoint Lie group,
and $H = H_{1} \times H_{2} \times \ldots \times H_{r}$ with $H_{i}$
simple. Let $\rho = \bigotimes_{i} \rho_{i}$ where $\rho_{i}$
is a finite dimensional irreducible complex representation of $H_{i}$.

\medskip

Assume that $\rho$ is even dimensional orthogonal. Then all $\rho$-liftable
automorphisms of $H$ are $\rho$-even if and only if one of the following case
happens.

\medskip

\textnormal{Case (1)}: Exactly one $\rho_{i}$ is even dimensional.
In this case all $\rho_{i}$-liftable automorphisms
of $H_{i}$ are $\rho_{i}$-even.

\medskip

\textnormal{Case (2)}: Exactly two, say $\rho_{i}$ and $\rho_{j}$,
are even dimensional, and moreover, either $\rho_{i}$ and $\rho_{j}$
are not quasi-equivalent, or ${\rm dim} (\rho_{i}) \equiv 0 \pmod{4}$.

\medskip

\textnormal{Case (3)}: At least three $\rho_{i}$s are even dimensional.
\end{theorem}

\medskip

\begin{corollary} \label{T:313}
All notations and assumptions are same as in Theorem ~\ref{T:312}. Then
$\rho$ is LFMO-special, if and only if each $\rho_{i}$ has weight $1$,
and one of the following holds:

\medskip

\textnormal{Case (1)}: Exactly one $\rho_{i}$ is even dimensional.
In this case $\rho_{i}$ is LFMO-special.

\medskip

\textnormal{Case (2)}: Same as \textnormal{Case (2)} in Theorem ~\ref{T:312}.

\medskip

\textnormal{Case (3)}: Same as \textnormal{Case (2)} in Theorem ~\ref{T:312}.
\end{corollary}

\medskip

\emph{Proof of Corollary of Corollary ~\ref{T:313} using Theorem ~\ref{T:312}}:
Note that if Case (1) occurs, then $\rho_{i}$ must be orthogonal since
al other $\rho_{j}$ are odd dimensional orthogonal, and $\rho$ is also orthogonal.
Now this corollary follows directly from Theorem ~\ref{T:312},
Proposition ~\ref{T:306} and Lemma ~\ref{T:307}.

\medskip

\qedsymbol

\medskip

\emph{Proof of theorem ~\ref{T:312}}:

\medskip

From Theorem ~\ref{T:308}, it suffices for us to
check the $\rho$-liftable automorphisms of $H$ in Type 1L and Type 2L.

\medskip

Step 1: Only if part. If exactly one $\rho_{i}$ is even dimensional,
and some $\phi_{i} \in {\rm Aut} (H_{i})$ is $\rho_{i}$-odd, then
$\phi' = 1 \times \ldots \times \rho_{i} \times \ldots \times 1$ is also
$\rho$ odd (Proposition ~\ref{T:311}, Type 1L); If exactly two, say
$\rho_{i}$ and $\rho_{j}$ are even dimensional, and $\rho_{i}$
and $\rho_{j}$ are quasi-equivalent, then there is a $\lambda: H_{i}
   \overset{\cong}{\to} H_{j}$ such that $T_{\lambda}$ is $\rho$-odd.
So the ``only if'' part follows.

\medskip

Step 2: Assume that all conditions in three cases hold. Want to
prove that all $\rho$-liftable automorphisms of $H$ in
Type 1L and Type 2L are $\rho$-even.

\medskip

Case (1) $+$ Type 1L: Consider
$\phi = \phi_{1} \times \ldots \times \phi_{r} \in {\rm Aut} (H)$
in as in Type 1L of Proposition ~\ref{T:311}. As $\phi_{i}$ is $\rho_{i}$-even
from the assumption of Case (1),
and all other $\phi_{j}$ are $\rho_{j}$-neutral as $\rho_{j}$
are odd dimensional, $\phi$ is then $\rho$-even (Proposition ~\ref{T:311}, Type 1L).

\medskip

Case (1) $+$ Type 2L: Consider
$T_{\lambda}$ as in Type 2L of Proposition ~\ref{T:311}
where $\lambda: H_{i} \overset{\cong}{\to} H_{j}$.
Then $\rho_{i}$ and $\rho_{j}$ are odd dimensional, and thus
for some unique $k \ne i, j$, $\rho_{k}$ is even dimensional.
Then $T_{\lambda}$ is $\rho$-even. (Proposition ~\ref{T:311}, Type 2L)

\medskip

Case (2) \& (3) $+$ Type 1L: Consider
$\phi = \phi_{1} \times \ldots \times \phi_{r} \in {\rm Aut} (H)$
in as in Type 1L of Proposition ~\ref{T:311}. As $\phi_{i}$ is $\rho_{i}$-even
and $\rho_{j}$ is $\rho_{j}$-even for $j \ne i$ from the assumption of Case (2)
and (3), $\phi$ is then $\rho$-even (Proposition ~\ref{T:311}, Type 1L).

\medskip

Case (2) $+$ Type 2L: Consider
$T_{\lambda}$ as in Type 2L of Proposition ~\ref{T:311}.
If $\rho_{i}$ and $\rho_{j}$ are even dimensional,
then, as ${\rm dim} (\rho_{i}) \equiv 0 \pmod{4}$,
$T_{\lambda}$ is $\rho$-even (Proposition ~\ref{T:311}, Type 2L).
If $\rho_{i}$ and $\rho_{j}$
are odd dimensional, then there must be two indices $k, l$
with $\rho_{k}$ and $\rho_{l}$ are even dimensional, and
thus $T_{\lambda}$ is also $\rho$-even (Proposition ~\ref{T:311}, Type 2L).

\medskip

Case (3) $+$ Type 2L: Consider
$T_{\lambda}$ as in Type 2L of Proposition ~\ref{T:311}.
No mater whether $\rho_{i}$ and $\rho_{j}$ are even dimensional,
there must be a $k \ne i, j$
with $\rho_{k}$ is even dimensional, and
thus $T_{\lambda}$ is also $\rho$-even (Proposition ~\ref{T:311}, Type 2L).

\medskip

Hence all $\rho$-liftable automorphisms of $H$ of Type 1L
and Type 2L are $\rho$-even. By Theorem ~\ref{T:308},
all $\rho$-liftable automorphisms are $\rho$-even.
Hence the theorem.

\qedsymbol

\medskip

\begin{proposition} \label{T:313Y}
Let $H$ be a connected complex reductive group and $\rho$ be a finite dimensional
complex representation. If $\rho (H)$ has trivial center, then
$\rho$ has weight $1$. The converse is true if $\rho$ is irreducible.
\end{proposition}

\medskip

\emph{Proof}: The weight set of $\rho$ forms a saturated set in $\Lambda$ (Propsotion
21.3 of \cite{Hu-GTM9}), the weight lattice
of $T$. Then the first assertion for $H$ semisimple follows from Lemma 13.4B of \cite{Hu-GTM9}.
For reductive $H$, since $\rho$ factors through its quotient which is semisimple of
adjoint type, then the first assertion for general case follows also from the semisimple case.

\medskip

For the converse, since $\rho$ is irreducible, $Z (H)$, the center of $H$
acts as scalars (Schur's Lemma, cf.\ Lemma ~\ref{T:206}), and hence we have the converse.

\qedsymbol

\medskip

Now we come to our main theorems.

\medskip

\emph{Proof of Theorem ~\ref{TM:D}}:

\medskip

We also proceed in steps. Recall that $H = H_{1} \times \ldots \times H_{r}$,
and $\rho = \bigotimes_{i} \rho_{i}$ where $\rho_{i}$ is a self-dual
representation of $H_{i}$.

\medskip

Step (1): $H$ is simply connected (resp.\ adjoint). i.e.,
$T = 1$ and all $H_{i}$ are simply connected (resp.\ adjoint).
Note that the center of $\rho (H)$ is trivial if and only if
$\rho$ has weight $1$ (Proposition ~\ref{T:313Y}).
Hence Theorem ~\ref{TM:D} is directly
from Corollary ~\ref{T:313}.

\medskip

Step (2): General semisimple $H = H_{1} \times \ldots \times H_{r}$.

\medskip

Let $\tilde{H} = \tilde{H}_{1} \times \ldots \tilde{H}_{r}$ where
$\tilde{H}_{i}$ is the simply connected cover of $H_{i}$, together
with the covering maps $\pi_{H_{i}}: \tilde{H}_{i} \to H_{i}$,
and $\pi_{H}: \tilde{H} \to H$ compatible to $\pi_{H_{i}}$.
Let $\tilde{\rho}_{i} = \rho_{i} \circ \pi_{H_{i}}$
and $\tilde{\rho} = \rho \circ \pi_{H} = \bigotimes_{i} \tilde{\rho}_{i}$.

\medskip

Then as $\pi_{H}$ is a finite central isogenies,
from the definition, we see easily that $\rho$ is LFMO-special,
if and only if $\tilde{\rho}$ is LFMO-special. In fact, all the three
conditions in the definitions are the same for $\rho$ and $\tilde{\rho}$.
Also, similarly the same for $\rho_{i}$ and $\tilde{\rho}_{i}$.

\medskip

Hence for each statement of the assertions
in all three cases in Theorem ~\ref{TM:D} are equivalent
for $H$ and for $\tilde{H}$. (Note that by our definition, $\rho_{i}$
and $\tilde{\rho}_{i}$ are quasi-equivalent!
So the assertion above works for Case (2).)
Then finally, Theorem ~\ref{TM:D} for $H$ is equivalent to Theorem ~\ref{TM:D}
for $\tilde{H}$, which is already done in Step (1).

\medskip

Step (3): General $H = T \times H_{0}$ where
$H_{0} = H_{1} \times \ldots \times H_{r}$ the derived group of $H$.
Now Theorem ~\ref{TM:D} follows from Step (2)
since $\rho$ is LFMO-special if and only $\rho |_{H_{0}}$ is so.

\qedsymbol

\medskip

\begin{theorem}  \label{T:314}
Let $H$ be a simple group of adjoint type and $\rho$ be the adjoint representation
of $H$. Then $\rho$ is LFMO-special if and only if
$H$ is of the type $A_{4 n} (n \ge 1), B_{2 n} (n \ge 1), C_{2 n} (n \ge 2)$
, $E_{6}, E_{8}, F_{4}, G_{2}$.
\end{theorem}

\medskip

This is also Theorem A in \cite{Wang2007}. See Table ~\ref{Table:1}.

\medskip

\emph{Proof}: Let $\mathfrak{h}$ be the Lie algebra of $H$, $\kappa$
the Killing form of $\mathfrak{h}$, $G = SO (\mathfrak(h), \kappa)$,
and $\Omega = O (\mathfrak{h}, \kappa)$. Then $\rho$ is orthogonal
and its image in $GL (\mathfrak{h})$ lies in $G$. Moreover, $\rho$
definitely has weight $1$. Thus, $\rho$ is LFMO-special
if and only if all $\phi \in {\rm Aut} (H)$ is $\rho$-even. In this case,
${\rm dim} (\rho) = {\rm dim} (\mathfrak{H})$ is even, and equivalently,
$H$ has even rank. So we focus on adjoint simple groups $H$ of type
$A_{2 n} (n \ge 1), B_{2 n} (n \ge 1), C_{2 n} (n \ge 2), D_{2 n} (n \ge 3)$,
$E_{6}, E_{8}, F_{4}, G_{2}$.

\medskip

Let $T$ be a maximal torus of $H$ with its Lie subalgebra
$\mathfrak{t} \in \mathfrak{h}$, $\Phi$ the root set of $\mathfrak{h}$
with respect to $\mathfrak{t}$ and $\Delta$ a simple basis.
Fix $\Sigma = (\mathfrak{t}, \{u_{\alpha}\} (\alpha \in \Delta))$
where $u_{\alpha}$ an eigenvector in $\mathfrak{h}$ of the root $\alpha$.

\medskip

Then by basic Lie theory,
\[
{\rm Aut} (H) = {\rm Int} (H) \rtimes {\rm Aut} (H, \Sigma)
\]
where ${\rm Aut} (H, \Sigma)$ is the set
of $\phi \in {\rm Aut} (H)$ that fixes $\Sigma$.
In particular, can choose an eigenbasis of $T$:
$(t_{\alpha} (\alpha \in \Delta), u_{\beta} (\beta \in \Phi))$
such that ${\rm Aut} (H, \Sigma)$ is a group
of permutations of coordinates. In particular,
${\rm Aut} (H, \Sigma) \in \Omega$.

\medskip

Since each $\phi \in {\rm Aut} (H)$ is $\rho$-liftable, and lifts
to $\phi$ itself when identify ${\rm Aut} (H)$ with ${\rm Aut} (\mathfrak{h})$
via the adjoint representation $\rho$, by the discussion above, we have,
to prove that all $\phi$ are $\rho$-even, it suffices
to prove that all $\phi \in {\rm Aut} (H, \Sigma)$
are even coordinate permutations.

\medskip

For each $\phi \in {\rm Aut} (H, \Sigma)$, $\phi$ permute
the positive root set $\Phi^{+}$ and the negative root set
$\Phi^{-} = - \Phi^{+}$ in exactly the same style. Hence
the signature of $\phi$ agrees with the signature of $\phi$
on $\Delta$, or equivalently, the Dynkin diagram.

\medskip

Hence for $H$ of even rank, $\rho$ is LFMO-special
if and only if all automorphisms
of $H$ are $\rho$-even, if and only if
all automorphisms in ${\rm Aut} (H, \Sigma)$ are even coordinate
permutations, if and only if all automorphisms of the
Dynkin diagram of $H$ are even permutations.
Now Table ~\ref{Table:1} will finally conclude our theorem.

\medskip

% Here inserts the Dynkin diagram.

%%\begin{picture}(50, 10)
%% Solid circles
%% \multiput(0, 0)(10, 0){5}{\circle*{2}}
%% Lines
%% \multiput(0, 0)(10, 0){3}{\line(1, 0){10}}
%% \multiput(30, -.5)(0, 1){2}{\line(1, 0){10}}
%% Label
%% \put(-3, -4){1}
%% \put(7, -4){2}
%% \put(17, -4){3}
%% \put(27, -4){4}
%% \put(37, -4){5}
%% Arrow
%% \put(34.5, -.8){$\rangle$}
%% \end{picture}

\begin{table}
\caption{Dynkin Diagrams of Even Rank}
\label{Table:1}
\begin{center}
\begin{tabular}{|c|c|c|c|}
\hline
{\tiny Type} & Diagram &
\tabincell{c}{\tiny Generators of \\ \tiny Automorphism \\ \tiny Group}
& \tiny{Signatures} \\ \hline
% A_{2 n}
\tabincell{c}{\tiny $A_{2 n}$ \\ \tiny $(n \ge 1)$} &
\begin{picture}(80, 10)(-5, -3)
% Node
\multiput(0, 0)(10, 0){3}{\circle*{2}}
\multiput(50, 0)(10, 0){3}{\circle*{2}}
% Line
\multiput(0, 0)(10, 0){3}{\line(1, 0){10}}
\multiput(40, 0)(10, 0){3}{\line(1, 0){10}}
% Dotted Line
\multiput(31, 0)(1, 0){9}{\circle*{0.5}}
% Label
\put(-3, -4){\tiny $1$}
\put(7, -4){\tiny $2$}
\put(17, -4){\tiny $2$}
\put(45, -4){\tiny $2n\!-\!2$}
\put(56, -4){\tiny $2n\!-\!1$}
\put(67, -4){\tiny $2n$}
\end{picture}
& \tabincell{c}{\tiny $(1 \ 2n) (2\ 2n-1)$
     \\ \tiny $ \cdots (n+1 \ n)$}
& \tabincell{c}{\tiny ${(-1)}^{n}$ \\ \tiny $\surd$ for $2 | n$
     \\ \tiny $\times$ for $2 \nmid n$}
\\[.5ex] \hline
% B_{2 n}
\tabincell{c}{\tiny $B_{2 n}$ \\ \tiny $(n \ge 2)$} &
\begin{picture}(80, 10)(-5, -3)
% Node
\multiput(0, 0)(10, 0){3}{\circle*{2}}
\multiput(50, 0)(10, 0){3}{\circle*{2}}
% Line
\multiput(0, 0)(10, 0){3}{\line(1, 0){10}}
\multiput(40, 0)(10, 0){2}{\line(1, 0){10}}
\multiput(60, -.3)(0, .6){2}{\line(1, 0){10}}
% Dotted Line
\multiput(31, 0)(1, 0){9}{\circle*{0.5}}
% Label
\put(-3, -4){\tiny $1$}
\put(7, -4){\tiny $2$}
\put(17, -4){\tiny $3$}
\put(45, -4){\tiny $2n\!-\!2$}
\put(56, -4){\tiny $2n\!-\!1$}
\put(67, -4){\tiny $2n$}
% Angle
\put(64.5, -.8){$\rangle$}
\end{picture}
& 1
& $\surd$
\\[.5ex] \hline
% C_{2 n}
\tabincell{c}{\tiny $C_{2 n}$ \\ \tiny $(n \ge 2)$} &
\begin{picture}(80, 10)(-5, -3)
% Node
\multiput(0, 0)(10, 0){3}{\circle*{2}}
\multiput(50, 0)(10, 0){3}{\circle*{2}}
% Line
\multiput(0, 0)(10, 0){3}{\line(1, 0){10}}
\multiput(40, 0)(10, 0){2}{\line(1, 0){10}}
\multiput(60, -.3)(0, .6){2}{\line(1, 0){10}}
% Dotted Line
\multiput(31, 0)(1, 0){9}{\circle*{0.5}}
% Label
\put(-3, -4){\tiny $1$}
\put(7, -4){\tiny $2$}
\put(17, -4){\tiny $3$}
\put(45, -4){\tiny $2n\!-\!2$}
\put(56, -4){\tiny $2n\!-\!1$}
\put(67, -4){\tiny $2n$}
% Angle
\put(64.3, -.8){$\langle$}
\end{picture}
& 1
& $\surd$
\\[.5ex] \hline
% D_{2 n}
\tabincell{c}{\tiny $D_{2 n}$ \\ \tiny $(n \ge 3)$} &
\begin{picture}(80, 20)(-5, -8)
% Node
\multiput(0, 0)(10, 0){3}{\circle*{2}}
\multiput(50, 0)(10, 0){2}{\circle*{2}}
\multiput(70, -5)(0, 10){2}{\circle*{2}}
% Line
\multiput(0, 0)(10, 0){3}{\line(1, 0){10}}
\multiput(40, 0)(10, 0){2}{\line(1, 0){10}}
\put(60, 0){\line(2, 1){10}}
\put(60, 0){\line(2, -1){10}}
% Dotted Line
\multiput(31, 0)(1, 0){9}{\circle*{0.5}}
% Label
\put(-3, -4){\tiny $1$}
\put(7, -4){\tiny $2$}
\put(17, -4){\tiny $3$}
\put(44, -4){\tiny $2n\!-\!3$}
\put(55, -4){\tiny $2n\!-\!2$}
\put(67, -9){\tiny $2n\!-\!1$}
\put(67, 1){\tiny $2n$}
\end{picture}
& \tabincell{c}{\tiny $(2n\!-\!1\ 2n)$ if $n \ne 4$ \\
\tiny $(3 \ 4),\,(1 \ 3 \ 4)$ if $n = 4$}
& \tabincell{c}{\tiny $-1$ if $n \ne 4$ \\ \tiny $-1$, $1$ if $n = 4$ \\
\tiny  $\times$ Always!}
\\[.5ex] \hline
% E_{6}
$E_{6}$ &
\begin{picture}(60, 20)(-5, -14)
% Node
\multiput(0, 0)(10, 0){5}{\circle*{2}}
\put(20, -10){\circle*{2}}
% Line
\multiput(0, 0)(10, 0){4}{\line(1, 0){10}}
\put(20, 0){\line(0, -1){10}}
% Label
\put(-3, 2){\tiny $1$}
\put(7, 2){\tiny $3$}
\put(17, 2){\tiny $4$}
\put(27, 2){\tiny $5$}
\put(37, 2){\tiny $6$}
\put(17, -14){\tiny $2$}
\end{picture}
& $(1 \ 6) (3 \ 5)$
& $1 \quad \surd$
\\[.5ex] \hline
% E_{8}
$E_{8}$ &
\begin{picture}(80, 20)(-5, -14)
% Node
\multiput(0, 0)(10, 0){7}{\circle*{2}}
\put(20, -10){\circle*{2}}
% Line
\multiput(0, 0)(10, 0){6}{\line(1, 0){10}}
\put(20, 0){\line(0, -1){10}}
% Label
\put(-3, 2){\tiny $1$}
\put(7, 2){\tiny $3$}
\put(17, 2){\tiny $4$}
\put(27, 2){\tiny $5$}
\put(37, 2){\tiny $6$}
\put(47, 2){\tiny $7$}
\put(57, 2){\tiny $8$}
\put(17, -14){\tiny $2$}
\end{picture}
& 1
& $\surd$
\\[.5ex] \hline
% F_{4}
$F_{4}$ &
\begin{picture}(50, 10)(-5, -3)
% Node
\multiput(0, 0)(10, 0){4}{\circle*{2}}
% Line
\multiput(0, 0)(20, 0){2}{\line(1, 0){10}}
\multiput(10, -.5)(0, 1){2}{\line(1, 0){10}}
% Label
\put(-3, -4){\tiny $1$}
\put(7, -4){\tiny $2$}
\put(17, -4){\tiny $3$}
\put(27, -4){\tiny $4$}
% Angle
\put(14.5, -.8){$\rangle$}
\end{picture}
& 1
& $\surd$
\\[.5ex] \hline
% G_{2}
$G_{2}$ &
\begin{picture}(40, 10)(-5, -3)
% Node
\multiput(0, 0)(10, 0){2}{\circle*{2}}
% Line
\multiput(0, -.6)(0, .6){3}{\line(1, 0){10}}
% Label
\put(-3, -4){\tiny $1$}
\put(7, -4){\tiny $2$}
% Angle
\put(4.5, -.8){$\langle$}
\end{picture}
& 1
& $\surd$
\\[.5ex] \hline
\end{tabular}
\end{center}
\end{table}

\bigskip

Note that finally we get $A_{4 n} (n \ge 1), B_{2 n} (n \ge 1), C_{2 n} (n \ge 2)$
, $E_{6}, E_{8}, F_{4}, G_{2}$.

\medskip

\qedsymbol

\medskip

\emph{Proof of Corollary ~\ref{TM:E}}: From the facts below about the simple
factor $H_{j}$ which is adjoint, we see
that the corollary follows directly from Theorem ~\ref{TM:D}
and Theorem ~\ref{T:314}.
First, adjoint representation of $H_{i}$
has weight $1$. Second, the adjoint representations of $H_{i}$
and $H_{j}$ are quasi-equivalent if and only if $H_{i}$ and $H_{j}$
are isomorphic. Third, ${\rm dim} (\rho_{j}) = {\rm dim} (H_{j})$
is even if and only if $H_{j}$ has even rank.

\medskip

\qedsymbol

\bigskip

\section{Odds and Ends} \label{S:4}

\subsection{$E_{8}$} \label{SS:401}

\ssmedskip

At the end of this notes, we quote the following theorem, which follows from
A result of Seitz (\cite{Seitz87}, \cite{Seitz91}).

\medskip

\begin{theorem} \label{T:401}
Let $H = \mPGL (3, \bC)$ and $G = E_{8} (\bC)$. Then there exists an embedding
$\rho: H \hookrightarrow G$ such that $m' (\rho; G) > 1$.
\end{theorem}

\qedsymbol

\emph{Remark}: There exists a homomorphism $i: {\rm Spin} (16, \bC) \hookrightarrow
    E_{8} (\bC)$ with kernel $\bZ / 2 \bZ$.
Let $\rho'_{1}: \mPGL (3, \bC) \to \mSL (8, \bC)$ be the adjoint representation, and
$\rho' = \rho'_{1} \oplus \rho'_{1}$. Then $\rho'$ is orthogonal type. Since
$Z (\mSL (3, \bC)) \cong \bZ / 3 \bZ$, $\rho'$ also lifts to $\tilde{\rho}':
    \mPGL (3, \bC) \to {\rm Spin} (16, \bC)$. Then $\rho = i \circ \tilde{\rho}'$.
From ~\ref{TM:C}, $\rho'$ is LFMO-special and
hence $m' (\tilde{\rho}; {\rm Spin} (16, \bC)) = 2$ (Theorem ~\ref{T:208},
Theorem ~\ref{T:214}).

\medskip

\subsection{Final Remark} \label{SS:402}

\ssmedskip

Since all arguments and statements are worked through complex linear algebraic groups,
and we don't involve any special topology, then everything and finally all results work
also for $K$, an algebraic closed field of characteristic $0$, for example,
$\bar{Q}_{l}$. Moreover, when $H$ and $G$ are connected, we have also Lie algebra version
of our concepts and results. In fact, some of well known results were first introduced
as Lie algebra version (see \cite{Dynkin}, \cite{Dynkin2}). Moreover,
we have also the compact real reductive group version. Our paper just focuses
on the version of the algebraic group over $\bC$.

\end{document}